%% file: SICON_Pathirajaetal_FINAL.tex
\documentclass[article,onefignum,onetabnum]{siamart190516}

\input{ex_shared}

\ifpdf
\hypersetup{
  pdftitle={Mckeanvlasov},
  pdfauthor={Sahani Pathiraja, Sebastian Reich, and Wilhelm Stannat}
}
\fi

\begin{document}

\maketitle

\begin{abstract}
  Various particle filters have been proposed over the last couple of decades with the common feature that the update step is governed by a type of control law.  This feature makes them an attractive alternative to traditional sequential Monte Carlo which scales poorly with the state dimension due to weight degeneracy.  This article proposes a unifying framework that allows us to systematically derive the McKean-Vlasov representations of these filters for the discrete time and continuous time observation case, taking inspiration from the smooth approximation of the data considered in  \cite{sr:crisan10} and \cite{Clark2005}.  We consider three filters that have been proposed in the literature and use this framework to derive It\^{o} representations of their limiting forms as the approximation parameter $\delta \rightarrow 0$.  All filters require the solution of a Poisson equation defined on $\mathbb{R}^{d}$, for which existence and uniqueness of solutions can be a nontrivial issue.  We additionally establish conditions on the signal-observation system that ensures well-posedness of the weighted Poisson equation arising in one of the filters.
\end{abstract}



	\section{Introduction}
\label{sec:intro}

Given a state process $\mathcal{X}_s$ which is indirectly observed through a second stochastic process $Z_s$, the goal of the standard filtering problem is to estimate the conditional distribution of the state $\mathcal{X}_t$ given observations until time $t$.  In the case where these processes evolve according to an Ornstein-Uhlenbeck process, the Kalman-Bucy filter (or Kalman filter for the discrete time setting) provides the exact posterior.  In the nonlinear setting, the Kushner-Stratonovich equation describes the evolution of the posterior distribution in continuous time, although it is rarely adopted in practice due to the intractibility of solving the stochastic partial differential equation (SPDE).  Consequently, a plethora of approximate algorithms have been proposed in the literature, most notably, Monte Carlo based methods that produce samples from the posterior.  For instance, ensemble Kalman filters (EnKFs) for state (e.g.,~\cite{sr:evensen03}, \cite{sr:tippett03}) and state-parameter estimation (e.g. \cite{Nusken2019}) have become popular as an approximate Bayesian computation style \cite{Nott2011} Monte Carlo based extension of the Kalman filter.  EnKFs have desirable accuracy and stability properties even for small ensemble sizes  (combined with techniques such as localisation and inflation \cite{Tong2016}, \cite{DelMoral2021}, \cite{Bishop2018}, \cite{deWiljes2020}), although it is inconsistent with Bayes' theorem in the nonlinear setting even as $N \rightarrow \infty$.  On the other hand, particle filters or sequential Monte Carlo algorithms which rely on importance sampling are theoretically consistent with the optimal filter as the number of particles $N \rightarrow  \infty$, although they have found limited use in practice due to well-known issues of weight degeneracy and the tendency to scale poorly with the dimension of the state space.  Over the last decade, there has been interest in developing Monte Carlo based algorithms which are consistent as $N \rightarrow \infty$ in the nonlinear setting but avoid such issues by steering particles according to a control law rather than relying on importance weights (see, for instance, \cite{Surace2019HowTA} in the continuous time setting).  Similar to \cite{Taghvaei2019}, we adopt the term ``particle flow'' filters, originally used by \cite{daumhunghomot} to refer to homotopy-based discrete time filters, to describe the type of aforementioned filter designed for the continuous or discrete time setting.  Such filters include, but are not limited to, the work of \cite{Mitter2003},\cite{Yang2011} in the fully continuous time setting and \cite{daumhuangparticleflow}, \cite{sr:reich10}, \cite{sr:crisan10}, \cite{Yang2014} for the case of discrete time observations.  The continuous signal-discrete observation case is handled either via the introduction of a homotopy to transport particles from prior to posterior in pseudotime (e.g.,~\cite{daumhunghomot}, \cite{sr:reich10}, \cite{Yang2014}) or through the construction of a single continuous time process by adopting a smooth approximation of the continuous time observation \cite{sr:crisan10}. \cite{Reich2019} also provides a mean-field equation for a particle flow filter in the fully continuous time setting for the case of smooth observations.  The main focus of this article is on the work of \cite{sr:crisan10} (hereafter referred to as the Crisan \& Xiong filter), the filter of \cite{sr:reich10} (hereafter referred to as the Reich filter) and the feedback particle filter (FPF) \cite{Yang2014}, \cite{Laugesen2015}.  We study the mean-field equations in each case, rather than the associated Monte Carlo approximations.

Throughout this article, we consider the following nonlinear filtering problem in continuous time 
\begin{align}
\label{eq:sigcts}
\text{signal process:} \quad d\mathcal{X}_t &= \mathcal{M}(\mathcal{X}_t)dt + dV_t,  \\
\label{eq:obscts}
\text{observation process:} \quad dZ_t &= h(\mathcal{X}_t)dt + dW_t, 
\end{align}
where $\mathcal{M}: \mathbb{R}^{d} \rightarrow \mathbb{R}^{d}$ and $h: \mathbb{R}^{d} \rightarrow \mathbb{R}$ with $\mathcal{M},h$ satisfying appropriate regularity conditions (e.g.,~globally Lipschitz), $\mathcal{X}_t \in \mathbb{R}^{d}$, $V_t$ and $W_t$ independent standard Wiener processes of dimension $d$ and $1$, respectively, and the filtration $\mathcal{Z}_t := \sigma(Z_s : s \leq t)$.  We focus on the scalar observation case without loss of generality, since the filters and our results extend to the multivariate observation case in a componentwise fashion.  Unless otherwise stated, all probability distributions throughout the article are absolutely continuous with respect to the Lebesgue measure and have finite second moment.  The filtering density, i.e.,~the conditional density of $\mathcal{X}_t$ given $\mathcal{Z}_t$, is denoted by $\theta_t$. 

The discrete time observation setting is more often considered in practice, particularly as continuous time measurements are rarely available.  Let $\{0 = t_1 < t_2 < \cdots < t_N = T\}$ be a partition of the time interval $[0,T]$ with fixed increment $\delta$, i.e.~$t_{n+1} - t_n = \delta$ for all  $n \in \{0, 1, \cdots, N\}$.   A typical discrete time observation model is of the form 
\begin{align}
	\label{eq:discobs}
	Y_n = h(\mathcal{X}_{t_n}) +  R^{1/2}\tilde{W}_n,
\end{align}
where $\tilde{W}_n$ is a standard Gaussian random variable independent of the signal noise, $R$ is the observation error variance, and $\mathcal{X}_{t_n}$ denotes the solution of \cref{eq:sigcts} at time $t_n$.  $Y_n$  can be viewed as an Euler-type approximation of the SDE \cref{eq:obscts} with step size $\delta$, i.e., 
\begin{align}
	\label{eq:ynmodel}
	\frac{\tilde{Z}_{n+1} - \tilde{Z}_n}{\delta } = h(\mathcal{X}_{t_n}) + \frac{1}{\sqrt{\delta}}  \tilde{W}_n, 
\end{align}
where $\tilde{Z}_n$ indicates the Euler approximation of $Z$ at the $n$th time step and $R = \frac{1}{\delta}$.  We consider an alternate perspective on the discrete time observation setting, as discussed further in \cref{sec:filtdisc}.  

\textbf{Crisan \& Xiong filter \cite{sr:crisan10}.} The Crisan \& Xiong filter  involves making a smooth approximation to the observing process \cref{eq:obscts} through a piecewise linear interpolation, which forms an alternative to the Euler-type approximation for continuous-discrete time filtering \cref{eq:discobs}-\cref{eq:ynmodel}.  It was developed as an extension of Clark's robust representation formula \cite{Clark2005} which ensures the filtering distribution depends continuously on the observation process (see also \cite{Coghi2021} for further extensions in the context of the ensemble Kalman filter).  Unlike the continuous-discrete FPF and Reich filter, the evolution of the filtering distribution occurs on a single time scale and is given by $\mu_t^\delta$, the conditional law of $X_t^\delta$ given $\mathcal{Z}_t$, which under conditions guaranteeing absolute continuity with respect to the Lebesgue measure has a density function $\rho_t^\delta (x)$.  The random variable $X_t^\delta$ evolves according to 
\begin{align}
	\label{eq:crisanmckeanvlasov}
	dX_t^\delta = \mathcal{M}(X_t^\delta)dt + dV_t + \frac{1}{\rho^\delta_t} \nabla u_t(X_t^\delta,\mu_t^\delta)dt,
\end{align}
where $\delta$ corresponds to the mesh size of the partition of the time interval $[0, T]$ used for generating the smooth approximation to the observation and 
\begin{align}
\label{eq:crisanfundsoln}
\nabla u_t(x,\mu) & = \frac{1}{\omega_{d}} \int \frac{(y-x)}{|x-y|^{d}} \left( m_t(y) - \overline{m}_t \right) \mu(dy), \\
\label{eq:crisanmt}
m_t(y) & :=  h(y)\frac{Z_{t_{n+1}} - Z_{t_n}}{\delta} - \frac{1}{2}h(y)^2 \quad \forall \enskip t\in [t_n, t_{n+1}),
\end{align}
where $\omega_{d}$ is the surface area of the unit ball in $\mathbb{R}^{d}$ and $\overline{m}_t = \int m_t(y) \mu(dy)$.  It is not difficult to see that \cref{eq:crisanfundsoln} corresponds to the fundamental solution (in gradient form) of the following Poisson equation defined on $\mathbb{R}^{d}$: 
\begin{align}
\label{eq:poisscrisan}
\nabla \cdot (\nabla u_t) = -(m_t - \overline{m}_t)\rho,
\end{align}
where $\mu(dx) = \rho(x)dx$ and $\rho$ is a smooth probability density function.  \cite{sr:crisan10} shows that $\rho_t^\delta$ is consistent with the solution of the Kushner-Stratonovich equation as $\delta \rightarrow 0$ under certain conditons on $\mathcal{M}$, $h$ and the initial distribution.  Additionally they establish the existence and uniqueness of such a process under further conditions, in particular, under the assumption of global Lipschitz continuity of $\frac{1}{\rho^\delta_t} \nabla u_t(x, \mu_t^\delta)$.  However, a precise limiting form of the stochastic process for $\delta \rightarrow 0$ (i.e.,~the continuous time observation limit) appears to be lacking.    

\textbf{Reich filter \cite{sr:reich10}.}  Here the continuous-discrete time filtering problem of the form of \cref{eq:discobs} is considered by constructing a deterministic transport.  The continuous time limit of a Bayes recursion with fixed observation $Y_n$ with $R = \frac{1}{\delta}$ is reformulated as a continuity equation, so that the transport from prior $\rho^-_n$ to posterior $\rho_n = \theta_n$ is governed by the ODE 
\begin{align}
\label{eq:reichodeorig}
\frac{dS_\lambda}{d\lambda} = M^{-1} \nabla \Pi_\lambda,
\end{align}  
evolving in pseudotime $\lambda \in [0,\delta]$.  $S_0$ is a random variable with associated density  $\rho_n^-$, and $\Pi_\lambda$ is the solution of the following Poisson equation: 
\begin{align}
\label{eq:reichpoiss}
\nabla \cdot (\rho_{n,\lambda} M^{-1} \nabla \Pi_\lambda) &= -(L_n - \overline{L}_{n,\lambda})\rho_{n,\lambda}, \\
\label{eq:reichlike} 
L_n(x) &:= h(x)Y_n - \frac{1}{2}h(x)^2,
\end{align}
where $\overline{L}_{n,\lambda} := \int L_n(x) \rho_{n,\lambda}(x) dx$, $\rho_{n,\lambda}$ indicates the probability density of the random variable $S_\lambda$ for observation sampling instance $n$ and $M \in \mathbb{R}^{d \times d}$ is a positive definite mass matrix that may potentially depend on $\lambda$.  As explained in \cite{sr:reich10}, the particular form of the right hand side of the ODE \cref{eq:reichodeorig} is the minimiser of the kinetic energy given by 
\begin{align*}
\mathcal{T}(v) = \frac{1}{2} \int v^T M v \rho_{n,\lambda} dx
\end{align*} 
and $\rho_{n,\lambda}(x) > 0$ for all $ x \in \mathbb{R}^{d}$.  Note the similarity between $L_n$ and $m_t$ as defined in \cref{eq:crisanmt}.  The choice of the mass matrix $M$ will be strongly application dependent, although from an approximation point of view, it can be seen as a kind of preconditioner for numerical solution of the elliptic equation (\ref{eq:reichpoiss}).   This approach is closely related to the homotopy-based formulations of \cite{Yang2014} and \cite{daumhunghomot}.     

\textbf{FPF {\cite{Yang2013}, \cite{Yang2014},\cite{Laugesen2015}}.} The FPF for continuous time observations is given by  
\begin{align}
\begin{aligned}
\label{eq:FPFctsmckean}
dX_t = \mathcal{M}(X_t)dt + dV_t + & K_t(X_t) \left( dZ_t - \frac{1}{2} \left( h(X_t) + \overline{h}_t \right) dt \right) \\
& + \frac{1}{2} \nabla K_t(X_t) K_t(X_t) dt,
\end{aligned}
\end{align}
where $\rho_t(x)$ indicates the conditional density of $X_t$ given $\mathcal{Z}_t$, $\overline{h}_t := \int h(x) \rho_t(x) dx$ and $K_t(x) \in \mathbb{R}^{d \times 1}$ with $K_t(x) := \nabla \phi_t(x)$, where $\phi_t(x)$ is given by the solution of the following weighted Poisson equation which must be solved for all $t$:
\begin{align}
\label{eq:poisscts}
\nabla \cdot (\rho_t(x) \nabla \phi_t(x)) = -(h(x) - \overline{h}_t ) \rho_t(x)
\end{align}
with the centring condition $\int \phi_t(x) \rho_t(x) dx = 0$ and where $\nabla K_t(X_t)$ indicates the Jacobian of $K_t$ evaluated at $x = X_t$.  
A derivation of the continuous time FPF via a variational formulation of Bayes' theorem was given in \cite{Laugesen2015} for the case of trivial signal dynamics.  Additionally, they show that the conditional density of $X_t$ given $\mathcal{Z}_t$ corresponds exactly to the solution of the Kushner-Stratonovich equation, $\theta_t$, under mild assumptions on the initial density and for bounded $h$.  

\cite{Yang2014} provides a formulation of the FPF for the setting where the signal process is given by \cref{eq:sigcts} and observations are available in discrete time according to \cref{eq:discobs}, i.e.,~the so-called continuous-discrete case. The signal dynamics and posterior calculation are split into two different processes; in particular the observations are incorporated through a homotopy.  Using the notation in \cite{Yang2014},  the signal process is first evolved from time $t_{n-1}$ to time $t_n$, with the final distribution having density denoted by $\rho^-_{n}$.  A homotopy is then constructed to evolve $\rho^-_{n}$ to the filtering density $\rho_n = \theta_{n}$ in pseudotime $\lambda$, i.e., 
\begin{align}
\label{eq:partflowS}
\frac{dS_\lambda}{d\lambda} = K_{\lambda} (S_\lambda) \left( Y_n - \frac{1}{2} \left(h(S_\lambda)+ \overline{h}_{n,\lambda} \right)  \right) +  \frac{1}{2} \Omega_\lambda (S_\lambda),
\end{align}
where $S_0$ is a random variable with associated density  $\rho_n^-$, $\lambda \in [0,\delta]$, $R = \frac{1}{\delta}$ in \cref{eq:discobs} and $\overline{h}_{n,\lambda} = \int h(x) \rho_{n,\lambda}(x) dx$.  Furthermore, $K_\lambda := \nabla \phi_\lambda$, where $\phi_\lambda$ is the solution of  
\begin{align}
\label{eq:ctsdiscpoiss1}
\nabla \cdot (\rho_{n,\lambda}(x) \nabla \phi_\lambda(x)) = -(h(x) - \overline{h}_{n,\lambda}) \rho_{n,\lambda}(x),
\end{align}
satisfying $\int \phi_\lambda(x) \rho_{n,\lambda}(x) dx = 0$ and $\Omega_\lambda := \nabla \xi_\lambda$, where $\xi_\lambda$ is the solution of 
\begin{align}
\label{eq:ctsdiscpoiss2}
\nabla \cdot (\rho_{n,\lambda}(x) \nabla \xi_\lambda(x)) = (r_\lambda(x) - \overline{r}_{n,\lambda}) \rho_{n,\lambda}(x), 
\end{align}
satisfying $ \int \xi_\lambda(x) \rho_{n,\lambda}(x) dx = 0$, where $r_\lambda:= \nabla \phi_\lambda \cdot \nabla h$, $\overline{r}_{n,\lambda} = \overline{h^2_{n,\lambda}} - (\bar{h}_{n,\lambda})^2$ and $\rho_{n,\lambda}(x)$ indicates the density of $S_\lambda$.  We note a minor sign difference between \cref{eq:ctsdiscpoiss2} and that presented in \cite{Yang2014}, which we justify through the derivation in \cref{sec:filtdisc} which is also applicable in the homotopy-based formulation.    The multivariate observation case is handled in a similar manner as for the continuous time case. 

A common feature of the aforementioned filters is the solution of a Poisson equation, generally defined on $\mathbb{R}^{d}$.  It is therefore important to establish conditions that ensure the existence and uniqueness of solutions in a particular function class.  It is well-known that the Poisson equation 
\begin{align}
\label{eq:poissstand}
\Delta u =- f
\end{align}
on an open and bounded domain $\Omega \subset \mathbb{R}^{d}$ for $f \in C(\Omega)$ has a unique solution for Dirichlet, Neumann or mixed boundary conditions  (in the standard filtering problem, the boundary conditions are of Neumann type).  Classical theory also provides conditions on $f$ guaranteeing solutions to \cref{eq:poissstand} on $\mathbb{R}^{d}$ and uniqueness in a function class depending on the regularity of $f$.  For weighted Poisson equations of the type $\nabla \cdot (\rho \nabla \phi) = f(\rho)$ arising in the FPF and Reich filters, it is typical to establish the existence and uniqueness of weak solutions in the $\rho$-weighted Sobolev space $H_\rho^1(\mathbb{R}^{d})$ combined with a centring condition.  It is clear that if the domain is bounded and with appropriate boundary conditions, one only needs that the weight function $\rho$ is smooth and strictly positive (to ensure the domain is connected).  Here we focus on weighted spaces with functions supported on $\mathbb{R}^{d}$ where the weight function satisfies a Poincar\'{e}  inequality.  Such a Poincar\'{e} inequality implies that the probability density has at least exponentially decaying tails \cite{Bakrygentled}, thereby supplying the crucial ingredient for solvability.  \cite{Laugesen2015} established well-posedness of \cref{eq:poisscts} via a Poincar\'{e} inequality for the case where $\mathcal{M} = 0$ and $h$ bounded with support $\mathbb{R}^{d}$.   

\subsection{Main Contributions}
As noted in both \cite{Taghvaei2019} and \cite{sr:reich10}, many particle flow filters require the solution of a type of first order PDE: 
\begin{align}
\label{eq:firstorderpdeK}
\nabla \cdot (\rho_t \mathcal{K}(x, \rho_t)) = F(x, \rho_t(x)),
\end{align}    
which has infinitely many solutions $\mathcal{K}$.  Our aim is to establish a unifying framework that allows us to systematically derive mean-field representations of various continuous time particle flow filters that are consistent with the optimal filter, only by varying the a priori assumptions on the form of $\mathcal{K}$.  Such a framework allows us to map continuously from the filtering equations in the discrete time observation to continuous time observation setting (see \cref{lem:meanfieldrep}), which is, to the best of our knowledge, missing in the literature.  The basis of this framework relies heavily on the smooth approximation of the data employed in \cite{sr:crisan10} and \cite{Clark2005}.  Additionally, we obtain representations of the limiting form of the Crisan \& Xiong  and Reich filters for the continuous time observation case, which has not yet been investigated in the literature.  The framework also allows us to derive an analogous form of the continuous-discrete FPF and Reich filter (hereafter referred to as the $\delta$-FPF and $\delta$-Reich filters, respectively) without the need for a two-step predict-update procedure.  However, we only provide a representation of these filters and conditions ensuring the existence of such representations are still to be determined.  The second main contribution of this article is to establish sufficient conditions on the signal-observation system that ensure well-posedness of the Poisson equations arising in the FPF and Crisan \& Xiong filters.  

The article is structured as follows.  Section 2 is devoted to deriving the aforementioned particle flow filters in the discrete time observation setting by considering a piecewise smooth approximation of the continuous time observation path.  The corresponding limiting equations for the continuous time observation setting ($\delta \rightarrow 0$) are established in section 3.  Sections 4 and 5 summarise the main results regarding well-posedness of the resulting Poisson equations for continuous time and discrete time observations, respectively.  Proofs of the well-posedness results are provided in section 6.

\subsection{Notation}
The following notation is used throughout the article.
\begin{itemize}
	\item $\mathbb{R}_{>0} = \{x \in \mathbb{R} \enskip | \enskip x > 0  \}$.
	\item $\succcurlyeq$ indicates matrix inequality. 
	\item $L^2_\rho(\mathbb{R}^d)$ denotes the Hilbert space of functions on $\mathbb{R}^d$ that are square integrable with respect to $\rho$.
	\item $H^k_{\rho}(\mathbb{R}^{d})$ denotes the Hilbert space of functions whose first $k$ derivatives (defined in the weak or distributional sense) are in $L^2(\mathbb{R}^d, p)$.  
	\item $H^k_{\rho,0}(\mathbb{R}^{d})$ denotes the space of all $f \in H^k_{\rho}(\mathbb{R}^{d})$ with the additional condition that $\int f(x) p(x) dx = 0$. 
	\item $C_b^k(\mathbb{R}^d)$ indicates the space of uniformly bounded $k$-times continuously differentiable functions with uniformly bounded derivatives up to $k$. 
	\item $\mathcal{P}$ indicates the space of all probability densities with finite first and second moments. 
	\item $\nabla$ and $\nabla^2$ indicate the gradient/Jacobian operator and Hessian operator, respectively, with respect to $x$.  The gradient is always assumed to be a column vector, and the Jacobian of a vector field $F:\mathbb{R}^n \rightarrow \mathbb{R}^m$ is an $n \times m$ matrix.  $D_\rho F(x)$ indicates the partial Fr\'{e}chet derivative operator of $F$ with respect to $\rho$ evaluated at $x$, and $D_\rho F(x)\xi$ indicates the operator acting on the function $\xi$.  Recall that if $U,W$ are Banach spaces, the Fr\'{e}chet derivative of a function $f: U \rightarrow W$ at $u \in U$ is a linear mapping $L: U \rightarrow W$ such that 
	\begin{align*}
	f(u + \xi) - f(u) - L(u)\xi = o(\norm{\xi})
	\end{align*}	
	for all $\xi \in U$.  The argument $x$ in $D_\rho F(x)\xi$ is occassionally suppressed throughout for notational brevity, although the meaning should be clear from context. 
	\item The overbar notation indicates an expectation, i.e.,~$\overline{h}_t = \int h(x) \rho_t(x) dx$, where unless otherwise stated, the probability density $\rho_t$ against which $h(x)$ is integrated is clear from context. The subscript $t$ is included to emphases the time-dependence of this quantity. 
	\item We adopt the shorthand notation $f_t^\delta$ for $f(x,\rho_t^\delta)$, likewise $f_t$ for $f(x, \rho_t)$ for the remainder of the article, except when the arguments of the function are emphasised.    
\end{itemize}	

\section{Filtering with smooth approximation to observations}
\label{sec:filtdisc}

Our strategy of obtaining continuous time filtering algorithms involves constructing a smooth observation path, either from discrete time observations or as an approximation of continuous time observations from \cref{eq:obscts}.  Such an approach allows the signal process and observation update to be considered in a single mean-field process, similar to \cite{sr:crisan10}.  Consider the time discretisation with increment $\delta$ discussed in \cref{sec:intro} and a smooth approximation of the Brownian path that generated the observations $t \rightarrow Z_t(\omega)$, $\omega \in \Omega$ via a piecewise linear interpolation $Z_t^\delta (\omega)$:
\begin{align}
	\label{eq:smoothobspath}
	Z_t^\delta (\omega) = Z_{t_n} (\omega) + \frac{Z_{t_{n+1}}(\omega) - Z_{t_n}(\omega)}{\delta } (t - t_n) \quad \forall \enskip t \in [t_n, t_{n+1}).
\end{align}  
In the remainder of the article we will drop the $\omega$ for notational ease. 
Consider then the following McKean-Vlasov SDE:  
\begin{align}
\label{eq:genMckeanvlasov}
dX_t^\delta = \mathcal{M}(X_t^\delta)dt + dV_t + a(X_t^\delta, \rho_t^\delta) dt + K(X_t^\delta, \rho_t^\delta) dZ_t^\delta,
\end{align}
where $a,K: \mathbb{R}^{d} \times \mathcal{P}(\mathbb{R}^{d}) \rightarrow \mathbb{R}^{d}$ and $\rho_t^\delta$ indicates the conditional density of $X_t^\delta$ given $\mathcal{Z}_t = \sigma(Z_s: s \leq t)$, assuming absolute continuity with respect to the Lebesgue measure.  It is readily seen that the corresponding (stochastic) Fokker-Planck equation is given by  	
\begin{align}
\label{eq:fokkerplankmv}
d \rho_t^\delta = \left( \mathcal{L^*} \rho_t^\delta - \nabla \cdot (\rho_t^\delta a(x, \rho_t^\delta)) \right) dt - \nabla \cdot (\rho_t^\delta K(x, \rho_t^\delta)) dZ_t^\delta, 
\end{align}
where $\mathcal{L}$ corresponds to the infinitesimal generator of \cref{eq:sigcts} and $\mathcal{L}^*$ its adjoint.  The goal is to obtain expressions for the coefficients $a$ and $K$ such that \cref{eq:fokkerplankmv} converges to the Kushner-Stratonovich equation as $\delta \rightarrow 0$.  A smooth approximation of the Kushner-Stratonovich equation (under appropriate conditions) is given by
\begin{align}
\label{eq:kushnerapprox}
d \theta_t^\delta = \left( \mathcal{L^*} \theta_t^\delta - \frac{1 }{2} (h^2 - \overline{h_t^2}) \theta_t^\delta\right) dt + (h - \overline{h}_t) \theta_t^\delta  dZ_t^\delta, 
\end{align}  
see  \cite{sr:crisan10}, \cite{Hu2002}.  Direct comparison of \cref{eq:fokkerplankmv} and \cref{eq:kushnerapprox} indicates that requiring   
\begin{align}
\label{eq:conta}
\nabla \cdot (\rho_t^\delta a(x, \rho_t^\delta)) = \frac{1 }{2} (h(x)^2 - \overline{h_t^2}) \rho_t^\delta, \\
\label{eq:contb}
\nabla \cdot (\rho_t^\delta K(x, \rho_t^\delta)) = - (h(x) - \overline{h}_t) \rho_t^\delta
\end{align}
achieves the desired convergence, given appropriate conditions on $a$ and $K$.  This is explored in further detail in \cref{sec:filtcts}.  
	\begin{remark}
	An alternative to the discrete time observation model \cref{eq:smoothobspath} that also allows us to construct a single ODE combining evolution under the signal dynamics and observation update is via mollification \cite{sr:reich10}.  The observation update is incorporated at all discrete time points when an observation becomes available through a Dirac delta function.  The discontinuities induced by the delta function are avoided by mollification.  
\end{remark}

\begin{remark}
	It is possible to consider further variations to the structure of the McKean-Vlasov process \cref{eq:genMckeanvlasov}, whilst ensuring its conditional law still corresponds to the filtering density.  A notable example in the linear-Gaussian case is the stochastically perturbed ensemble Kalman--Bucy filter (EnKBF)  \cite{Bergemann2012}, which involves the addition of a noise term with the same properties as the observation noise to \cref{eq:genMckeanvlasov}.  A further stochastic extension of the EnKBF can be found in \cite{8618878}.  \cite{daumstochpart} also considers a stochastic extension which is valid for smooth strictly positive filtering densities.  A stochastic variant of the FPF is also presented in \cite{Reich2019}.
\end{remark}

\noindent It is clear that both $a$ and $K$ are not uniquely defined by \cref{eq:conta} and \cref{eq:contb}, respectively, which leads to various formulations of nonlinear filters.  We explore some of these in further detail below.

\textbf{Crisan \& Xiong filter \cite{sr:crisan10}.} By direct comparison of \cref{eq:crisanmckeanvlasov} and \cref{eq:genMckeanvlasov}, we have that 
\begin{align*}
\frac{1}{\rho^\delta_t} \nabla u_t^\delta dt =  a_t^\delta dt + K_t^\delta dZ_t^\delta,
\end{align*}
from which we obtain 
\begin{align}
\label{eq:crisident}
\nabla \cdot \nabla u_t^\delta dt = \nabla \cdot (\rho_t^\delta a_t^\delta)dt + \nabla \cdot (\rho_t^\delta K_t^\delta) dZ_t^\delta. 
\end{align}
Differentiating the smooth approximation \cref{eq:smoothobspath} gives $\frac{dZ_t^\delta}{dt} = \frac{Z_{t_{n+1}} - Z_{t_n}}{\delta}$, which together with \cref{eq:crisanmt}, \cref{eq:poisscrisan} and \cref{eq:crisident} implies that the coefficients are assumed to be of the form 
\begin{align*}
a_t^\delta = \frac{\nabla \alpha_t^\delta}{\rho_t^\delta} \enskip ; \quad  K_t^\delta = \frac{\nabla \beta_t^\delta}{\rho_t^\delta},
\end{align*}
where $\rho_t^\delta(x) > 0$ for all $x \in \mathbb{R}^{d}$ and $\alpha_t^\delta, \beta_t^\delta$ are the solutions of 
\begin{align}
\label{eq:pde_alpha_crisandisc}
\nabla \cdot (\nabla \alpha_t^\delta) &= \frac{1 }{2} (h^2 - \overline{h_t^2}) \rho_t^\delta, \\
\label{eq:pde_beta_crisandisc}
\nabla \cdot (\nabla \beta_t^\delta) &= -(h - \overline{h}_t)\rho_t^\delta. 
\end{align}
Notice also that this coincides with the required \cref{eq:conta} and \cref{eq:contb}.	
Regularity of the solutions $\alpha_t^\delta$ and $\beta_t^\delta$, as well as of the density $\rho_t^\delta$ are strongly related to the $\delta \rightarrow 0$ case and are discussed in further detail in \cref{sec:filtcts} and \cref{sec:wellposedcts}.  
 
\textbf{$\delta$-Reich filter.} A McKean-Vlasov formulation of the homotopy-based filter in \cite{sr:reich10} can be obtained by replacing the ODE model \cref{eq:reichodeorig} with 
\begin{align*}
dX_t^\delta = \mathcal{M}(X_t^\delta)dt + dV_t + M_t^{-1} \nabla \Pi_t^\delta dt, 
\end{align*}
where $M$ is a $\mathbb{R}^{d \times d}$ positive definite mass matrix.  Again, direct comparison of the above with \cref{eq:genMckeanvlasov} and the associated Fokker-Planck equations implies 
\begin{align*}
\nabla \cdot (\rho_t^\delta M_t^{-1} \nabla \Pi_t^\delta) dt = \nabla \cdot (\rho_t^\delta  a_t^\delta)dt + \nabla \cdot (\rho_t^\delta  K_t^\delta) \dot{Z}_t^\delta dt, 
\end{align*}
which together with \cref{eq:conta} and \cref{eq:contb} implies that the coefficients in \cref{eq:genMckeanvlasov} are assumed to be of the form 
\begin{align}
\label{eq:reichaK}
a_t^\delta = M_t^{-1} \nabla \Omega_t^\delta \enskip ; \quad K_t^\delta = M_t^{-1} \nabla \Lambda_t^\delta,
\end{align}
where $\Omega_t^\delta$ and $\Lambda_t^\delta$ are the solutions of 
\begin{align}
\label{eq:poissomeg}
\nabla \cdot (\rho_t^\delta M_t^{-1} \nabla \Omega_t^\delta) = \frac{1 }{2} (h^2 - \overline{h_t^2}) \rho_t^\delta, \\
\label{eq:poisslambd}
\nabla \cdot (\rho_t^\delta M_t^{-1} \nabla \Lambda_t^\delta) = - (h - \overline{h}_t) \rho_t^\delta, 
\end{align}
respectively. Note that the weighted Poisson equation arising from \cref{eq:poissomeg} and \cref{eq:poisslambd} for $\Pi_t^\delta = \Omega_t^\delta + \Lambda_t^\delta\dot{Z}_t^\delta$ has the same form as \cref{eq:reichpoiss}.  
 
\textbf{$\delta$-FPF.} The homotopy-based FPF for discrete observations \cite{Yang2014} can be reformulated similar to the $\delta$-Reich filter.  Specifically, the coefficients in \cref{eq:genMckeanvlasov} are assumed to be of the form 
\begin{align}
\label{eq:fpfcoeffdisc}
a_t^\delta = -\frac{1}{2} K_t^\delta(h + \overline{h}_t) + \frac{1}{2} \nabla \psi_t^\delta \enskip ; \quad K_t^\delta = \nabla \phi_t^\delta.
\end{align}
The expression for $\phi_t^\delta$ follows directly from \cref{eq:contb}:
\begin{align}
\label{eq:poissphidisc}
\nabla \cdot (\rho_t^\delta \nabla \phi_t^\delta) = - (h - \overline{h}_t) \rho_t^\delta. 
\end{align}
An expression for $\psi_t^\delta$ is obtained by rewriting \cref{eq:conta} as
\begin{subequations}
	\begin{align}
	\nabla \cdot ( {\rho}_{t}^\delta a_t^\delta ) &= \frac{\rho_t^\delta}{2} \left(  (h - \overline{h}_t)(h + \overline{h}_t) + (\overline{h}_t)^2 - \overline{h^2_t} \right) \\
	& = -\frac{1}{2}  (h + \overline{h}_t) \nabla \cdot \left(K_t^\delta \rho_t^\delta \right)  + \frac{\rho_t^\delta}{2}((\overline{h}_t)^2 - \overline{h^2_t} ) \\
	\label{eq:pdeamanip}
	& = - \frac{1}{2} \nabla \cdot \left(K_t^\delta (h + \overline{h}_t) \rho_t^\delta \right) + \frac{\rho_t^\delta}{2} \nabla h^T K_t^\delta + \frac{\rho_t^\delta}{2}((\overline{h}_t)^2 - \overline{h^2_t} ), 
	\end{align}
\end{subequations}
combined with the assumed form of $a_t^\delta$ gives 
\begin{align}
\label{eq:poisspsidisc}
\nabla \cdot (\rho_t^\delta \nabla \psi_t^\delta) &= \left( \nabla h^T K_t^\delta + (\overline{h}_t)^2 - \overline{h^2_t}  \right) \rho_t^\delta. \\
&= (r_t - \bar{r}_t) \rho_t^\delta 
\end{align} 
where $r_t:= \nabla h^T K_t^\delta$.  The last equality can be seen by taking the weak form of (\ref{eq:poissphidisc}) with $h$ as the test function, which then gives $\bar{r}_t:=\int \nabla h^T K_t^\delta \rho_t^\delta dx = \overline{h^2_t} -  (\overline{h}_t)^2$.  Combining all leads to    
\begin{align}
\begin{aligned}
\label{eq:mehtadisc}
dX_t^\delta = \mathcal{M}(X_t^\delta)dt + dV_t &+  \nabla \phi(X_t^\delta, \rho_t^\delta) \left( dZ_t^\delta -\frac{1}{2} (h(X_t^\delta) + \overline{h}_t  )dt \right) \\
&+ \frac{1}{2} \nabla \psi(X_t^\delta, \rho_t^\delta) dt, 
\end{aligned}
\end{align}
where $\psi(X_t^\delta, \rho_t^\delta)$ is the solution of \cref{eq:poisspsidisc}.  Conditions ensuring well-posedness of \cref{eq:poissphidisc} and \cref{eq:poisspsidisc} are established in \cref{sec:wellposedctsdisc}.  
\begin{remark}
	\label{rem:correcctsdisc}
	Despite the modification to the formulation of the filter, the underlying structure of the Poisson equations is unchanged (cf. \cref{eq:poissphidisc}  with \cref{eq:ctsdiscpoiss1} and \cref{eq:poisspsidisc} with \cref{eq:ctsdiscpoiss2}).  
\end{remark} 

 	\section{Filtering with continuous time observations}
 \label{sec:filtcts}
 
 We first state an important assumption regarding the continuous time filtering distribution at time $t$, which ensures that it admits a density with respect to the Lebesgue measure, and that it is the unique solution in $W^{k,2}(\mathbb{R}^{d})$ of \cref{eq:kushnerstratfokkerplanck} (see, e.g.,~Theorem 7.11 and 7.17 in \cite{sr:crisan}).  
 
 \begin{assumption}
 	\label{ass:condsolnKS}
 	The following conditions on the signal-observation system \cref{eq:sigcts}-\cref{eq:obscts} are satisfied:
 	\begin{itemize}
 		\item The initial density $\rho_0 \in W^{k,2}(\mathbb{R}^{d})$ with finite second moment.
 		\item The signal drift $\mathcal{M}: \mathbb{R}^{d} \rightarrow \mathbb{R}^{d}$ and observation drift $h: \mathbb{R}^{d} \rightarrow \mathbb{R}$ are $C_b^{k+1}$ functions. 
 	\end{itemize}
 \end{assumption}
 \noindent 
These conditions are also sufficient to guarantee the existence of a unique solution to \cref{eq:kushnerapprox}, as proved in \cite{Hu2002}.  The starting point for the analysis in this section is the following McKean-Vlasov SDE in Stratonovich form: 
 \begin{align}
 \label{eq:meanfieldstrat}
 dX_t = \mathcal{M}(X_t)dt + dV_t + a(X_t, \rho_t)dt + K(X_t, \rho_t) \circ dZ_t, 
 \end{align}
 with 
 \begin{align}
 \label{eq:conteqa}
 \nabla \cdot (\rho_t a_t) &= \frac{1 }{2} (h^2 - \overline{h_t^2}) \rho_t, \\
 \label{eq:conteqK}
 \nabla \cdot (\rho_t K_t) &= - (h - \overline{h}_t) \rho_t. 
 \end{align}
 where $V_t$ is a standard Brownian motion and $\rho_t$ is the conditional density of $X_t$ given $\mathcal{Z}_t$.  The $\circ$ notation in \cref{eq:meanfieldstrat} indicates the Stratonovich interpretation which must be considered for both arguments of $K$, i.e.,~$x$ and $\rho$.  \cref{lem:kushnerrho} shows that the evolution equation for $\rho_t$ coincides with the Kushner-Stratonovich equation under the following assumption.
 \begin{assumption}
 	\label{eq:lebesguerho}
 	The conditional law of $X_t$ given $\mathcal{Z}_t$, where $X_t$ satisfies \cref{eq:meanfieldstrat}, is absolutely continuous with respect to the Lebesgue measure with density given by $\rho_t$.
 \end{assumption}
 \begin{lemma}
 	\label{lem:kushnerrho}
 	\textbf{Evolution equation for $\rho_t$.}  For the system \cref{eq:meanfieldstrat} and \cref{eq:conteqa}-\cref{eq:conteqK}, and given \cref{eq:lebesguerho}, $\rho_t$ evolves according to the Kushner-Stratonovich equation, i.e.
 	\begin{align}
 	\label{eq:kushnerrhot}
 	d\rho_t = \mathcal{L}^*\rho_tdt + (h - \overline{h}_t)\rho_t (dZ_t - \overline{h}_tdt).
 	\end{align}
 	
 	\begin{proof}
 		The nonlinear Fokker-Planck equation for the conditional density of $X_t$ given $\mathcal{Z}_t$, $\rho_t$, where $X_t$ is the solution of \cref{eq:meanfieldstrat} is given by 
 		\begin{subequations}
 			\begin{align}
 			d \rho_t &= \left( \mathcal{L}^*\rho_t - \nabla \cdot (\rho_t a(x, \rho_t)) \right) dt - \nabla \cdot (\rho_t K(x, \rho_t)) \circ dZ_t  \\
 			\label{eq:fokkerplanckmckeanstrat}
 			& = \left( \mathcal{L}^*\rho_t -\frac{1}{2}(h^2 - \overline{h_t^2}) \rho_t \right) dt + (h - \overline{h}_t)\rho_t \circ dZ_t. 
 			\end{align}
 		\end{subequations}
 		Using the well-known It\^{o} -Stratonovich conversion between any two semimartingales we have 
 		\begin{align*}
 		\int_0^t (h - \overline{h}_s) \rho_s \circ dZ_s = \int_0^t (h - \overline{h}_s) \rho_t dZ_s + \frac{1}{2} \int_0^t d((h - \overline{h}_s) \rho_s) dZ_s. 
 		\end{align*} 
 		Evaluating the second term on the r.h.s of the above using It\^{o}'s  product rule gives 
 		\begin{align*}
 		d((h - \overline{h}_s) \rho_s)dZ_s &= hd\rho_sdZ_s - \left( \rho_s d\overline{h}_s + \overline{h}_s d\rho_s + d\rho_s d\overline{h}_s   \right)dZ_s \\
 		& = \left( h(h - \overline{h}_s)\rho_s  - \overline{h}_s(h - \overline{h}_s)\rho_s \right) ds - \rho_s \int h d\rho_s dx dZ_s   \\
 		& =  (h - \overline{h}_s)^2\rho_s ds - \rho_s(\overline{h_s^2} - (\overline{h}_s)^2) ds \\
 		& = (h^2 - \overline{h_s^2})\rho_s ds - 2\overline{h}_s(h - \overline{h}_s) \rho_s ds. 
 		\end{align*}
 		Combining with \cref{eq:fokkerplanckmckeanstrat} gives 
 		\begin{subequations}
 			\begin{align}
 			d \rho_t &=  \mathcal{L}^*\rho_tdt -\frac{1}{2}(h^2 - \overline{h_t^2}) \rho_t dt + (h - \overline{h}_t)\rho_t dZ_t + \frac{1}{2} (h^2 - \overline{h_t^2})\rho_t dt - \overline{h}_t(h - \overline{h}_t) \rho_t dt \\
 			\label{eq:kushnerstratfokkerplanck}
 			& =  \mathcal{L}^*\rho_tdt + (h - \overline{h}_t)\rho_t (dZ_t - \overline{h}_tdt).
 			\end{align} 
 		\end{subequations}
 	\end{proof}
 \end{lemma}
 \noindent
 Our main result is \cref{lem:meanfieldrep}, which provides an It\^{o} representation of \cref{eq:meanfieldstrat}, for which we rely on the following assumption.  
 
 \begin{assumption}
 	\label{ass:frechetdiff}
 	Given a specific functional form of $K:\mathbb{R}^{d} \times L^1(\mathbb{R}^{d}) \rightarrow  \mathbb{R}^{d}$ satisfying \cref{eq:conteqK}, $K(x, \rho)$ is twice continuously Fr\'{e}chet differentiable in $x$ and $\rho$.  
 \end{assumption}
 
 \begin{lemma}
 	\label{lem:meanfieldrep}
 	\textbf{McKean-Vlasov It\^{o}  representation.} 
 	Under \cref{ass:frechetdiff}, \cref{eq:meanfieldstrat} has the It\^{o}  representation 
 	\begin{align}
	\begin{aligned}
	 \label{eq:meanfield}
 	dX_t  =  \mathcal{M}(X_t)dt + dV_t +& K(X_t, \rho_t)  \left(dZ_t -\frac{1}{2}(h(X_t) + \overline{h}_t)dt \right) \\
 	& + \frac{1}{2}(\nabla K(X_t, \rho_t) )^T K(X_t, \rho_t)dt + \mathcal{J}(X_t, \rho_t) dt, 
 	\end{aligned}
 	\end{align}
 	where 
 	\begin{align}
 	\label{eq:defnJ}
 	\mathcal{J}(x, \rho_t) &:=  a(x, \rho_t) + g(x, \rho_t)  +   \frac{1}{2}K(x, \rho_t) (h(x) + \overline{h}_t), \\
 	\label{eq:defng}
 	g(x, \rho_t) &:= \frac{1}{2} D_\rho K(x, \rho_t)\rho_t(h(x) - \overline{h}_t), 
 	\end{align}
 	and $a: \mathbb{R}^{d} \times L^1(\mathbb{R}^{d}) \rightarrow  \mathbb{R}^{d}$, $K:\mathbb{R}^{d} \times L^1(\mathbb{R}^{d}) \rightarrow  \mathbb{R}^{d}$ and $g: \mathbb{R}^{d} \times L^1(\mathbb{R}^{d}) \rightarrow  \mathbb{R}^{d}$, where $a(x, \rho_t)$ and $K(x, \rho_t)$ satisfy \cref{eq:conteqa} and \cref{eq:conteqK} respectively.  Furthermore, 
 	\begin{align}
 	\label{eq:pdeJ}
 	\nabla \cdot (\rho_t \mathcal{J}(x, \rho_t)) = 0.
 	\end{align}
 	
 	\begin{proof}
 		To obtain the It\^{o}  formulation of \cref{eq:meanfieldstrat}, we again rely on the It\^{o} -Stratonovich relation for two semimartingales, giving  
 		\begin{align*}
 		\int_0^t K(X_s, \rho_s) \circ dZ_s = \int_0^t K(X_s, \rho_s) dZ_s + \frac{1}{2} \int_0^t dK(X_s, \rho_s) dZ_s. 
 		\end{align*}
 		Let $K^{i}_t$ be the $i$-th component of the vector field $K_t$.  Since $K^i_t$ is a twice continuously Fr\'{e}chet differentiable scalar function and the driving Wiener process in \cref{eq:kushnerstratfokkerplanck} is finite dimensional, we can apply the standard It\^{o} formula to obtain  
 		\begin{align*}
 		dK^{i}_t = (\nabla K^{i}_t)^T dX_s + \frac{1}{2} \text{Tr} \left( K^T_t (\nabla^2 K^{i}_t) K_t \right)  ds & + D_\rho K^{i}_t d\rho_s \\
 		& + \frac{1}{2}D^2_\rho K_t^{i}\rho_s^2(h - \overline{h}_s)^2ds,
 		\end{align*}
 		where the notation $D_yFh$ and $D^2_yFh$ indicate the first and second order partial Fr\'{e}chet derivative operators respectively of $F$ w.r.t the variable $y$ acting on $h$.  Combining the above, together with \cref{eq:kushnerrhot} gives 
 		\begin{align*}
 		K(X_t, \rho_t) \circ dZ_t = K(X_t, \rho_t)  dZ_t &+ \frac{1}{2}(\nabla K(X_t, \rho_t) )^T K(X_t, \rho_t) dt \\
 		&+ \frac{1}{2}D_\rho K(X_t, \rho_t)\rho_t(h - \overline{h}_t)dt.
 		\end{align*}
 		Under \cref{ass:frechetdiff}, we obtain the Fr\'{e}chet derivative of $K(x, \rho_t(x))$ with respect to $\rho_t$ acting on the perturbation $\rho_t(h - \overline{h}_t)$ as		
 		\begin{align*}
 		\nabla \cdot \left(\rho_t D_\rho K(x, \rho_t)\rho_t(h - \overline{h}_t)\right) & = -\nabla \cdot (K \rho_t (h - \overline{h}_t)) + \rho_t \int h(h - \overline{h}_t) \rho_t dx - (h - \overline{h}_t)^2 \rho_t \\
 		& = -\nabla \cdot (K \rho_t) (h - \overline{h}_t) - \nabla h^T K \rho_t + \rho_t (\overline{h_t^2} - (\overline{h}_t)^2) - (h - \overline{h}_t)^2 \rho_t \\
 		& = - \left( \nabla h^T K + (\overline{h}_t)^2 -\overline{h_t^2} \right)\rho_t.
 		\end{align*}
 		Since $g(x, \rho_t):= \frac{1}{2} D_\rho K(x, \rho_t)\rho_t(h - \overline{h}_t)$, we have shown that $g$ satisfies 
 		\begin{align}
 		\label{eq:pdeg}
 		\nabla \cdot (\rho_t g_t) = - \frac{1}{2}\left( \nabla h^T K_t + (\overline{h}_t)^2 -\overline{h_t^2} \right)\rho_t. 
 		\end{align}
 		Recall that $a_t$ and $K_t$ satisfy \cref{eq:conteqa} and \cref{eq:conteqK} respectively.  It follows from a similar line of reasoning as in \cref{eq:pdeamanip} that $a_t + g_t$ satisfies 
 		\begin{align}
 		\label{eq:pdeag}
 		\nabla \cdot (\rho_t (a_t + g_t)) = -\frac{1}{2} \nabla \cdot (\rho_t K_t(h + \overline{h}_t) ).
 		\end{align}
 		The definition of $\mathcal{J}(x, \rho_t)$ \cref{eq:defnJ} combined with the above directly implies \cref{eq:pdeJ}. 
 	\end{proof}
 \end{lemma}
 
 Although there is a vast literature on convergence theorems for smooth approximations of SDEs and SPDEs, there has been relatively little investigation in the context of McKean-Vlasov SDEs.  It was shown in \cite{sr:crisan10} that when 1) $\rho_0$ has finite second moment and is strictly positive and 2) $\mathcal{M}$ and $h$ are at least $C_b^2$ functions, then the solution of \cref{eq:fokkerplankmv} at time $t$, $\rho_t^\delta$, converges to the solution of the Kushner Stratonovich equation as $\delta \rightarrow 0$.  Similar results were obtained in \cite{Hu2002} for the Zakai equation.  
 Given \cref{lem:kushnerrho}, it is not unreasonable to expect that under certain conditions, \cref{eq:genMckeanvlasov} converges to \cref{eq:meanfieldstrat}  as $\delta \rightarrow 0$, at least weakly.  We defer the convergence analysis and and identification of sufficient conditions to a future publication (see \cite{Pathiraja2020} for some preliminary work in the case of ordinary SDEs with unbounded coefficients).  Here we focus purely on providing a representation of the filters in the formal $\delta \rightarrow 0$ limit.  Specifically, given (\cref{ass:frechetdiff}), \cref{lem:meanfieldrep} can be used to obtain the limiting McKean-Vlasov representations in It\^{o} form of the filters discussed in \cref{sec:filtdisc}, as detailed below.  \\
 \textbf{Crisan \& Xiong filter:} 
 Recall the assumed forms on the coefficients 
 \begin{align}
 \label{eq:crisancoeffscont}
 a_t = \frac{\nabla \alpha_t}{\rho_t}; \quad K_t = \frac{\nabla \beta_t}{\rho_t},
 \end{align}
 where $\alpha_t$ and $\beta_t$ are the solutions of 
 \begin{align}
 \label{eq:pde_alpha_crisan}
 \nabla \cdot (\nabla \alpha_t) &= \frac{1 }{2} (h^2 - \overline{h_t^2}) \rho_t, \\
 \label{eq:pde_beta_crisan}
 \nabla \cdot (\nabla \beta_t) &= -(h - \overline{h}_t)\rho_t. 
 \end{align}
 There exists solutions $\alpha_t$, $\beta_t$ in $C^2 \cap L^\infty$ whenever $h, \rho_t \in C^2 \cap L^\infty$ at least for $d \geq 3$ (see for example \cite{sr:crisan10}), which is ensured under \cref{ass:condsolnKS} and additional conditions (see Theorem 7.12 in \cite{sr:crisan}).  See also \cref{sec:wellposedcts} for further discussion on this issue. Additionally, we require $\rho_t > 0$, which is guaranteed under \cref{ass:condsolnKS} and when $\rho_0 > 0$ also \cite{sr:crisan10}.
 In order to apply \cref{lem:meanfieldrep}, we assume continuous Fr\'{e}chet differentiability of $\beta$ with respect to $\rho$.  Evaluating \cref{eq:defng} yields 
 \begin{align*}
 g_t &= -\frac{1}{2}\frac{\nabla \beta_t}{\rho_t} (h - \overline{h}_t) + \frac{1}{2}\frac{1}{\rho_t} D_\rho \nabla \beta_t \rho_t(h - \overline{h}_t) \\
 & := -\frac{1}{2}\frac{\nabla \beta_t}{\rho_t} (h - \overline{h}_t) + \frac{1}{\rho_t} \nabla q_t, 
 \end{align*}
 where $q_t:= \frac{1}{2}D_\rho \beta_t\rho_t (h - \overline{h}_t)$.  When $\beta_t \in C^2$ and is continuously Fr\'{e}chet differentiable with respect to $\rho_t$, we have symmetry of mixed partial Fr\'{e}chet derivatives, thereby permitting the change of differentiation order in the above.  In this case $\mathcal{J}_t$ takes the form 
 \begin{align*}
 \mathcal{J}_t = \frac{1}{\rho_t} \left(\nabla \alpha_t + \nabla q_t + \overline{h}_t {\nabla \beta_t}    \right),
 \end{align*}  
 and \cref{eq:pdeJ} implies the following Laplace equation for $f_t:= q_t + \alpha_t + \overline{h}_t \beta_t$:
 \begin{align}
 \nabla \cdot \nabla f_t = 0.
 \end{align}
 Additionally, from \cref{eq:pdeg} we have that $q_t$ satisfies  
 \begin{align}
 \label{eq:pde_q_crisan}
 \nabla \cdot (\nabla q_t) &= -\frac{\rho_t}{2} ((h - \overline{h}_t)^2 + (\overline{h}_t)^2 -\overline{h_t^2} ). 
 \end{align}
 Since there exists a solution in $C^2 \cap L^\infty$ to \cref{eq:pde_q_crisan}, \cref{eq:pde_alpha_crisan} and \cref{eq:pde_beta_crisan} when $h \in C^2 \cap L^\infty$, we have that $f_t \in C^2 \cap L^\infty$.  It then follows from Liouville's theorem that $f_t$ is a constant, so that $\nabla f_t = 0$.  That is, when $h \in C^2 \cap L^\infty$, $\mathcal{J}_t = 0$ so that the limiting form of the Crisan \& Xiong filter is given by 
 \begin{align*}
 dX_t = \mathcal{M}(X_t)dt + dV_t  &+ K(X_t,\rho_t) \left( dZ_t  -\frac{1}{2} (h(X_t) + \overline{h}_t) dt \right) \\
 &+  \frac{1}{2} \nabla K(X_t,\rho_t)^T K(X_t,\rho_t) dt, 
 \end{align*}
 with $K_t = \frac{\nabla \beta_t}{\rho_t}$.  
 
 \textbf{Reich filter.}  Recall the  assumed forms of the coefficients 
 \begin{align}
 \label{eq:reichaKcont}
 a_t= M_t^{-1} \nabla \Omega_t \enskip ; \quad K_t = M_t^{-1} \nabla \Lambda_t,
 \end{align}
 where $\Omega_t$ and $\Lambda_t$ are the solutions of 
 \begin{align}
 \label{eq:pdereichomeg}
 \nabla \cdot (\rho_t M_t^{-1} \nabla \Omega_t) = \frac{1 }{2} (h^2 - \overline{h_t^2}) \rho_t, \\
 \label{eq:pdereichlam}
 \nabla \cdot (\rho_t M_t^{-1} \nabla \Lambda_t) = - (h - \overline{h}_t) \rho_t. 
 \end{align}
 Evaluating $g$ using \cref{eq:defng} and the above, and again assuming continuous Fr{\'e}chet differentiability gives  
 \begin{subequations}
 	\begin{align}
 	g_t&= \frac{1}{2} D_\rho K_{t}\rho_t(h - \overline{h}_t) \\
 	& = \frac{1}{2} D_\rho M_t^{-1}\rho_t(h - \overline{h}_t) \nabla \Lambda_t + \frac{1}{2}M_t^{-1} \nabla D_\rho \Lambda_t\rho_t (h - \overline{h}_t) \\
 	\label{eq:reichg}
 	&=:\frac{1}{2}\tilde{\Gamma}_t \nabla \Lambda_t + M_t^{-1} \nabla \tilde{\Theta}_t
 	\end{align}
 \end{subequations}
 where $\tilde{\Theta}_t := \frac{1}{2}D_\rho \Lambda_t\rho_t(h - \overline{h}_t)$ and $\tilde{\Gamma}_t := D_\rho M^{-1}_t\rho_t(h - \overline{h}_t)$.  Then making use of \cref{eq:reichaKcont}, $\mathcal{J}_t$ takes the form 
 \begin{align}
 \label{eq:Jreich}
 \mathcal{J}_t  = M_t^{-1} \nabla \left( \Omega_t + \tilde{\Theta}_t \right) + \frac{1}{2}  \left( \tilde{\Gamma}_t   +  (h + \overline{h}_t) M_t^{-1}  \right) \nabla \Lambda_t,
 \end{align} 
 
 \noindent We now consider three possibilities for $M_t$, i.e., 
 
 \textit{Case} $M_t = \rho_t I$. It follows trivially in this case that the filter is equivalent to the Crisan \& Xiong filter, since $a_t$ and $K_t$ then have the same structure as \cref{eq:crisancoeffscont} and the corresponding PDEs \cref{eq:pdereichomeg} and \cref{eq:pdereichlam} collapse down to \cref{eq:pde_alpha_crisan} and \cref{eq:pde_beta_crisan} respectively. 
 
 \textit{Case} $M_t= I$.  It is clear that in this case, $\tilde{\Gamma}_t = 0$, so that 
 \begin{align*}
 \mathcal{J}_t =  \nabla \left( \Omega_t + \tilde{\Theta}_t \right) + \frac{1}{2}   (h + \overline{h}_t) \nabla \Lambda_t, 
 \end{align*} 
 where 
 \begin{align}
 \label{eq:pdethetaomg}
 \nabla \cdot (\rho_t  \nabla (\Omega_t + \tilde{\Theta}_t))  = -\frac{1}{2} \nabla \cdot (\rho_t (h + \overline{h}_t) \nabla \Lambda_t  ). 
 \end{align}
 In this case, the limiting form is given by 
 \begin{align*}
 dX_t = \mathcal{M}(X_t)dt + dV_t  + K(X_t, \rho_t) dZ_t  &+  \frac{1}{2} \nabla K(X_t,\rho_t)^T K(X_t,\rho_t) dt \\
 & + \nabla \left( \Omega + \tilde{\Theta} \right)(X_t, \rho_t) dt, 
 \end{align*}
 where $K(x,\rho_t) = \nabla \Lambda(x,\rho_t)$.  We emphasise that \cref{eq:pdethetaomg}, i.e., 
 	\begin{align*}
 	\nabla \cdot (\rho_t  \underbrace{\left[\nabla (\Omega_t + \tilde{\Theta}_t)  + \frac{1}{2} (h + \overline{h}_t) \nabla \Lambda_t \right]}_{= \mathcal{J}_t} ) = 0,
 	\end{align*}
 	does not imply that $\mathcal{J}_t= 0$.  This is because the term $\frac{1}{2} (h + \overline{h}_t) \nabla \Lambda_t $ is not of gradient type, meaning that the initial assumptions on the form of the coefficients $a_t$ and $K_t$ \cref{eq:reichaKcont} are not compatible with this particular form.  This emphasises the need for the specific structure assumed for $a$ in the FPF i.e.,~\cref{eq:fpfcoeffdisc} and \cref{eq:fpfacont} in order to achieve $\mathcal{J}_t = 0$, when $M_t = I$ and $K_t$ is of gradient type. 
 
\textit{Case} $M_t = P_t^{-1}$.  In this case, we have that for any admissible perturbation $\xi$, 
 \begin{equation}
 \label{eq:frechetPt}
 \begin{aligned}
 D_\rho P_t\xi &= \int (x - \overline{x}_t) (x - \overline{x}_t)^T \xi dx - \int (x - \overline{x}_t) \left( \int x \xi dx  \right)^T \rho_t dx \\
 & - \int \left(\int x \xi dx \right) (x - \overline{x}_t)^T \rho_t dx \\
 & =  \int (x - \overline{x}_t) (x - \overline{x}_t)^T \xi dx, 
 \end{aligned}
 \end{equation}
 which then implies
 \begin{align*}
 \tilde{\Gamma}_t = \int  (x - \overline{x}_t) (x - \overline{x}_t)^T (h-\overline{h}_t)\rho_t dx. 
 \end{align*}
 The calculations for this case are also relevant for the discussion of the linear-Gaussian case in \cref{sec:lineargauss}.
 
 \textbf{FPF \cite{Laugesen2015}.}
 Recall the assumed form for the coefficients in the FPF: 
 \begin{align}
 \label{eq:fpfacont}
 a_t &= -\frac{1}{2} K_t (h + \overline{h}_t) + \frac{1}{2} \nabla \psi_t, \\
 \label{eq:fpfKcont}
 K_t&= \nabla \phi_t, 
 \end{align}
 with 
 \begin{align}
 \label{eq:poiss_psi_mehta}
 \nabla \cdot (\rho_t \nabla \psi_t) &=  \left( \nabla h^T K_t+ (\overline{h}_t)^2 - \overline{h_t^2} \right)\rho_t, \\
 \label{eq:poiss_phi_mehta}
 \nabla \cdot (\rho_t \nabla \phi_t) &= -(h - \overline{h}_t)\rho_t. 
 \end{align}
 If $\phi(x, \rho_t)$ is continuously Fr\'{e}chet differentiable, \cref{eq:defng} and \cref{eq:pdeg} give
 \begin{align}
 \nonumber
 g(x, \rho_t) &= \frac{1}{2} \nabla \tilde{q}(x, \rho_t), 
 \end{align}
 where $\tilde{q}_t := D_\rho \phi_t\rho_t (h - \overline{h}_t)$ and $\tilde{q}_t$ satisfies 
 \begin{align}
 \label{eq:poissfrechmeta}
 \nabla \cdot (\rho_t \nabla \tilde{q}_t) &= - \left( \nabla h^T K_t + (\overline{h}_t)^2 -\overline{h_t^2} \right)\rho_t.
 \end{align}
 Then $\mathcal{J}_t$ takes the form 
 \begin{align*}
 \mathcal{J}_t = \frac{1}{2} \left( \nabla \psi_t + \nabla \tilde{q}_t \right). 
 \end{align*}
 If there exists a unique solution to \cref{eq:poiss_psi_mehta} in $H_{\rho_t, 0}^1(\mathbb{R}^{d})$ (which is ensured for the system \cref{eq:sigprocess}-\cref{eq:obsprocess}, see \cref{sec:wellposedcts} for further conditions), then both \cref{eq:poiss_psi_mehta}, \cref{eq:poissfrechmeta} and the fact that $\tilde{q}_t \in H_{\rho_t, 0}^1(\mathbb{R}^{d})$ imply that $\mathcal{J}_t = 0$.  Therefore, the FPF takes the form   
 \begin{align}
 \begin{aligned}
 \label{eq:familiarFPF}
 dX_t =  \mathcal{M}(X_t)dt + dV_t & + K(X_t,\rho_t)\left( dZ_t -\frac{1}{2} (h(X_t) + \overline{h}_t) dt \right) \\
 & + \frac{1}{2}(\nabla K(X_t,\rho_t))^T K(X_t,\rho_t)dt,
 \end{aligned}
 \end{align}
 where $K(x,\rho_t) = \nabla \phi(x,\rho_t)$.  The calculations above emphasise the importance of interpreting the Stratonovich $\circ$ in (\cite{Yang2011} and \cite{Laugesen2015}) as being with respect to space only.  \cite{8100938} obtained a derivation of the FPF for a continuous time system where the hidden state evolves on a matrix Lie group according to a Stratonovich SDE.  It would be of interest to extend the above analysis for such a manifold setting.  \\
 \\
 There are several factors to consider when choosing the form of $a_t$ and $K_t$, such as well-posedness of the corresponding PDEs as well as of the resulting McKean-Vlasov SDE in the continuous time limit. For instance, $a_t$ and $K_t$ should satisfy certain regularity estimates to ensure the integrals in (\ref{eq:meanfield}) are well-defined, i.e.,
 \begin{align*}
 \int \norm{K_t}^2 \rho_tdx < \infty,
 \end{align*}	 
 which ensures that the stochastic integral $\int K(X_t, \rho_t) dZ_t$ is well-defined in $L^2$ and 
 \begin{align*}
 \int \norm{\nabla K_t^T K_t} \rho_t dx < \infty, \enskip \int \norm{a_t} \rho_t dx < \infty, \enskip \int \norm{\mathcal{J}_t} \rho_t dx < \infty
 \end{align*}	
 for the Lebesgue integrals.   Additionally, nongradient forms of $a_t,K_t$ can give rise to Poisson equations which have a unique classical solution under fairly mild conditions on $\mathcal{M}, h$, although establishing well-posedness of the McKean-Vlasov SDE may be nontrivial due to the $\frac{1}{\rho_t}$ term.  From a numerical/computational perspective, non-gradient forms of $K_t$ may require explicit estimates of the posterior density, unlike in the case of gradient form of $K_t$ where the density enters as an average and therefore does not require an explicit estimate when Monte Carlo is employed. These issues will be investigated in further detail in a future publication. 
 
 \section{Filter representations in the linear-Gaussian case}
 \label{sec:lineargauss}

 In this section we investigate the form of the filters in the linear-Gaussian setting: 
 \begin{align}
 	\label{eq:lineargausssig}
 	d\mathcal{X}_t = A\mathcal{X}_t + dV_t, \\
 	\label{eq:lineargaussobs}
 	dZ_t = H\mathcal{X}_t + dW_t,
 \end{align}
 where $A:\mathbb{R}^d \rightarrow \mathbb{R}^d$ and $H: \mathbb{R}^d \rightarrow \mathbb{R}$.
 
 \subsection{Smooth approximation to observations}
  When the initial density is Gaussian, the conditional density $\rho_t^\delta$ arising from the smooth approximation to the observation process (\ref{eq:lineargaussobs}) is Gaussian for all $t$.  In particular, the FPF and Reich filters coincide and take the same form as the EnKBF \cite{Bergemann2012}, i.e.,
 \begin{align}
 	\label{eq:enkbfformdisc}
 	dX_t^\delta = AX_t^\delta dt + dV_t + P_tH^T \left(dZ_t^\delta -\frac{1}{2}H(X_t^\delta + \overline{X}_t^\delta)dt \right), 
 \end{align}
 as demonstrated below and also in \cite{Yang2014} and \cite{sr:reich10} for the homotopy-based formulations.  Note also that the EnKBF also coincides with a constant gain approximation of the FPF \cite{Taghvaei2017}.  
 
 \textbf{Crisan \& Xiong filter \cite{sr:crisan10}.}  Using the fact that $\rho_t^\delta$ is Gaussian, it follows that $\beta_t = P_t H \Phi_t(x)$ is a solution of \cref{eq:pde_beta_crisandisc} for $d = 1$, where $\Phi_t(x)$ indicates the Gaussian cumulative distribution function, i.e.,~$\int_{-\infty}^x \rho_t^\delta(y) dy$.  It is then immediate that $K_t^\delta = \frac{\nabla \beta_t^\delta}{\rho_t^\delta} = P_tH$, where $P_t \in \mathbb{R}$ indicates the variance and $H: \mathbb{R} \rightarrow \mathbb{R}$.  This fact can be used to determine solutions to \cref{eq:pde_alpha_crisandisc}; following similar calculations as in the FPF leads to 
 \begin{align*}
 	\frac{d^2\alpha_t^\delta }{dx^2} = \frac{d}{dx} \left( H(x+ \overline{x}_t) \frac{d \beta_t^\delta}{dx}  \right),
 \end{align*}	
 from which it is clear that $\frac{d \alpha_t^\delta}{dx} = H(x + \overline{x}_t)P_tH \rho_t^\delta + C$, where $C = 0$ if we require $|\alpha_t^\delta| \rightarrow 0$ as $|x| \rightarrow \infty$.  This then leads to the familiar EnKBF form \cref{eq:enkbfformdisc}.  A similar expression does not appear to hold for $d \geq 2$ (for further discussion on this issue, see \cite{Yang2013}). 
 
 \textbf{$\delta$-Reich filter.} For the specific choice of the mass matrix $M_t^{-1} = P_t$, where for the remainder of the article $P_t$ indicates the covariance of the mean-field process at time $t$, it is clear that $ \Lambda_t^\delta = (x - \overline{x}_t)^TH^T$ is a solution of 
 \begin{align}
 	\label{eq:reichKlineargauss}
 	\nabla \cdot (\rho_t^\delta P_t \nabla \Lambda_t^\delta) = - H(x - \overline{x}_t) \rho_t^\delta. 
 \end{align}
 Furthermore, it follows that $\Omega_t^\delta = -\frac{1}{4} (x + \overline{x}_t)^TH^TH(x + \overline{x}_t)$ is a solution of 
 \begin{align*}
 	\nabla \cdot (\rho_t^\delta P_t \nabla \Omega_t^\delta) = \frac{1 }{2} H(x x^T - \overline{(x x^T)}_t)H^T \rho_t^\delta, 
 \end{align*}
 since 
 \begin{align*}
 	\nabla \cdot (\rho_t^\delta P_t \nabla \Omega_t^\delta) & = \rho_t^\delta \nabla \cdot (P_t \nabla \Omega_t^\delta) + \nabla \rho_t^\delta \cdot \nabla \Omega_t^\delta \\
 	& = -\frac{1}{2}\rho_t^\delta \nabla \cdot (P_t H^T H(x + \overline{x}) + \frac{1}{2} (P_t^{-1}(x - \overline{x}_t))^T P_t(H^T H(x + \overline{x}_t)) \rho_t^\delta \\
 	& = -\frac{1}{2} \rho_t^\delta Tr(P_t H^T H) + \frac{1}{2} (x - \overline{x}_t)^T H^TH(x + \overline{x}_t) \rho_t^\delta \\
 	& = \frac{1}{2} (x^TH^TH x - H\overline{(x x^T)}_t H^T)\rho_t^\delta \\
 	& = \frac{1 }{2} H(x x^T - \overline{(x x^T)}_t)H^T \rho_t^\delta. 
 \end{align*}
 Substituting in $a_t^\delta = M_t^{-1} \nabla \Omega_t^\delta, \enskip  K_t^\delta = M_t^{-1} \nabla \Lambda_t^\delta$ in \cref{eq:genMckeanvlasov} leads to \cref{eq:enkbfformdisc}.
 
 \textbf{$\delta$-FPF.} It follows from direct substitution that $\phi_t^\delta = (x - \overline{x}_t)^TP_t H^T $ is a solution of \cref{eq:poissphidisc}, which then implies that \cref{eq:poisspsidisc} takes the form
 \begin{align}
 	\label{eq:bvppsi}
 	\nabla \cdot (\rho_t^\delta \nabla \psi_t^\delta) = \rho_t^\delta H P_t H^T - \rho_t^\delta HP_tH^T = 0, 
 \end{align}
 for which $\psi_t^\delta = 0$ is a solution.  Furthermore, both $\phi_t^\delta = (x - \overline{x}_t)^TP_t H^T $ and $\psi_t^\delta = 0$ are unique solutions of their respective PDEs in the function class $H^1_{\rho_t^\delta, 0} (\mathbb{R}^{d})$, as detailed in \cref{sec:wellposedctsdisc}.  Substituting into \cref{eq:mehtadisc} gives \cref{eq:enkbfformdisc}.

 \subsection{Continuous time observations}
 Again, we assume the initial density $\rho_0$ is Gaussian.  It is well-known that the solution of the Kushner-Stratonovich equation is then a Gaussian density. The limiting form of the Reich filter coincides with the EnKBF, since with the specific choice $M_t^{-1} = P_t$, it can easily be seen that $ \Lambda_t = \Lambda_t^\delta$ and $\Omega_t = \Omega_t^\delta$.  The only difference to the $\delta > 0$ case is to determine the form of $\mathcal{J}_t$.  Letting $\xi = \rho_tH(x - \overline{x}_t)$ in \cref{eq:frechetPt} gives $\tilde{\Gamma}_t = 0$ 
 when $\rho_t$ is a multivariate Gaussian.  By a similar calculation as in \cref{eq:frechetPt}, we have that 
 \begin{align*}
 	\nabla \tilde{\Theta}_t & = \frac{1}{2} \nabla \left(\int x^T H(x - \overline{x}_t)H^T  \rho_t dx \right) \\
 	& = 0.
 \end{align*}
Combining these results with \cref{eq:Jreich} implies that $\mathcal{J}_t = 0$, 
 which then leads to the familiar EnKBF form. 
 The FPF in this case is known to be equivalent to the EnKBF (see, e.g.,~\cite{Yang2013}, \cite{Taghvaei2017}), which can also be seen by following the same arguments as for the $\delta > 0$ case and using the fact that the term $(D_xK)^TK = 0$ for $K_t = P_tH^T$.  Similarly, the limiting form of the Crisan \& Xiong filter coincides with the EnKBF for $d =1$ from the same reasoning as in the $\delta > 0$ by confirming that $\beta_t = P_tH^T \Phi_t(x)$ is a solution of \cref{eq:pde_beta_crisan}.  As mentioned earlier, we only provide the representations of these filters; the convergence of the smooth approximations with unbounded coefficients in this setting will be examined in future research.


 \section{Well-posedness of the Poisson equations for continuous time observations $\delta \rightarrow 0$}
 \label{sec:wellposedcts}
 
 We begin by focusing on the weighted Poisson equation arising in the FPF.  A crucial ingredient for establishing well-posedness here is the existence of a Poincar\'{e} inequality for the conditional density $\rho_t$ with domain = $\mathbb{R}^{d}$.  Since we have that $\rho_t  = \theta_t$ for all $t$ when $\rho_0 = \theta_0$ due to \cref{lem:kushnerrho}, we actually require a Poincar\'{e} inequality for the filtering density $\theta_t$.  Throughout this section, we will use $\rho_t$ in place of $\theta_t$ for consistency with the Poisson equations described earlier.  The following lemma details conditions guaranteeing a Poincar\'{e} inequality uniformly in time, which is of interest in its own right in addition to the solvability issue. \cite{VanHandelthesis2007} and \cite{Stannat2005} (see Sections 4.3.3 and 1.3, respectively) obtained related results for this system in the context of stability of the optimal filter.  Here we adopt a different approach for the proof and are able to obtain slightly weaker conditions on the signal-observation drift functions than in \cite{VanHandelthesis2007} (see \cref{rem:vhcond}).

 \begin{lemma}
 	\label{lem:PIctstime}
 	\textbf{Poincar\'{e}  inequality for continuous time case}.  Consider a probability space $(\Sigma, \mathcal{F}, P)$ with the usual conditions.  Given stochastic processes $\mathcal{X}_s, Z_s$ defined on this space and evolving according to
 	\begin{align}
 	\label{eq:sigprocess}
 	d\mathcal{X}_t = \nabla U(\mathcal{X}_t)dt + dV_t, \\
 	\label{eq:obsprocess}
 	dZ_t = H\mathcal{X}_tdt + dW_t, 
 	\end{align}
 	where $U: \mathbb{R}^{d} \rightarrow \mathbb{R}$, $H \in \mathbb{R}^{1 \times d}$, $U \in C^{4}$ and $dV_t$, $dW_t$ are independent Brownian increments and $\mathcal{X}_0 \sim \mu_0$ with density $\rho_0$ with respect to the Lebesgue measure.  Suppose the initial density is of the form $\rho_0(x) = \exp(-\mathcal{G}_0(x))$.  Consider $c_r, c_u, c_g \in \mathbb{R}_{>0}$. If the conditions  
 	\begin{enumerate}[label=(C\arabic*)]
 		\item \label{cond:a1}  $U$ is strongly concave with parameter $c_u$, i.e.~$\nabla^2 U \preccurlyeq -c_u I$  
 		\item \label{cond:a2}  $\check{\mathcal{G}}_0:= U + \mathcal{G}_0$, $\check{\mathcal{G}}_0: \mathbb{R}^d \rightarrow \mathbb{R}$ is strongly convex with parameter  $c_g$, i.e.~$\nabla^2 \check{\mathcal{G}}_0 \succcurlyeq c_g I$ 
 		\item \label{cond:a3}  $R(x):= \norm{Hx}^2 + (\Delta U + \norm{\nabla U}^2)(x) $, $R: \mathbb{R}^d \rightarrow \mathbb{R}$  is also strongly convex with parameter $c_r$  
 		\item \label{cond:a4} $\nabla U$ satisfies a linear growth condition, i.e.~there exists a $D < \infty$ such that $|\nabla U(x)| \leq D(1 + |x|)$ 
 	\end{enumerate}
 	are satisfied, then the conditional density $\rho_t(x)$ for the signal-observation pair satisfies a Poincar\'{e}  inequality for all $ t \geq 0$, i.e.
 	\begin{align}
 	\label{eq:PIlogconcont}
 	\int |f(x)|^2 \rho_t(x) dx \leq \kappa \int |\nabla f(x)|^2 \rho_t(x) dx, 
 	\end{align}    
 	for all test functions $f \in H_{\rho_t, 0}^1(\mathbb{R}^{d})$, and where $\kappa = \left(c_u + \min \left( c_g, \sqrt{\frac{c_r}{2}}  \right) \right)^{-1}$.  
 	
 \end{lemma} 

 The proof of the above lemma is given in \cref{sec:proofPIlemmacts}.  Conditions \ref{cond:a1}-\ref{cond:a3} are required to ensure log-concavity of the posterior density, whilst \ref{cond:a4} is needed to ensure the signal process possesses a unique strong solution and that change of measure via Girsanov holds (i.e.~to satisfy Novikov's condition).  The existence and uniqueness of weak solutions to the weighted Poisson equation \cref{eq:poiss_phi_mehta} also follows from classical methods once $\rho_t$ is known to satisfy a Poincar\'{e} inequality, as shown in the following theorem.  
 
 \begin{theorem}
 	\label{theo:wellposedcts}
 	\textbf{Well-posedness of Poisson equation in continuous time FPF.} Consider the signal-observation pair and conditions specified in \cref{lem:PIctstime}.  Then there exists a unique solution $\phi \in H_{\rho_t,0}^1(\mathbb{R}^{d})$ for all $t \geq 0$ that satisfies the following Poisson equation in weak form: 
 	\begin{align}
 	\label{eq:poissweakcts}
 	\int \nabla \phi(x) \cdot \nabla \psi(x) \rho_t(x) dx = \int H(x - \overline{x}_t) \psi(x) \rho_t(x) dx, 
 	\end{align}
 	for all test functions $\psi \in H_{\rho_t,0}^1(\mathbb{R}^{d})$. 
 	
 	\begin{proof}
 		The standard $H^1_{\rho_t}$ norm $\norm{u}_{H_{\rho_t}^1}:= \left( \int |u |^2 \rho_t(x) dx +  \int |\nabla u |^2 \rho_t(x) dx \right)^{1/2} $ and the norm $\norm{u}_{H_{\rho_t,0}^1}:= \left( \int |\nabla u |^2 \rho_t(x) dx \right)^{1/2}$ are equivalent, since by \cref{lem:PIctstime} 
 		\begin{align*}
 		\norm{u}_{H_{\rho_t,0}^1} \leq \norm{u}_{H_{\rho_t}^1} \leq c\norm{u}_{H_{\rho_t,0}^1},
 		\end{align*}
 		where $c > 0$. This implies that the l.h.s of \cref{eq:poissweakcts} is an inner product on ${H_{\rho_t,0}^1}$.  Furthermore, the r.h.s of \cref{eq:poissweakcts} is a bounded linear functional of $\psi(x)$ since 
 		\begin{align*}
 		|T(\psi)| &:= \left|\int (Hx - \overline{H}_t)\rho_t^{1/2} \psi(x)\rho_t^{1/2}dx \right| \\
 		& \leq   \left(\int (Hx-\overline{H}_t)^2 \rho_t dx\right)^{1/2} \left(\int |\psi|^2 \rho_t dx\right)^{1/2}  \\
 		& \leq C  \left(\int |\psi|^2 \rho_t dx\right)^{1/2}  \quad  (\rho_t \enskip \text{has finite first \& second moment)}\\
 		& \leq C_2  \left( \int |\nabla \psi(x)|^2 \rho_t dx \right)^{1/2}  \quad \text{(using \cref{lem:PIctstime})}\\
 		& = C_2 \norm{\psi(x)}_{H^1_{\rho_t,0}},
 		\end{align*}
 		Therefore, by the Riesz representation theorem, there exists a unique solution $\phi(x) \in H_{\rho_t,0}^1(\mathbb{R}^{d})$. 
 	\end{proof}
 	
 \end{theorem}
 
 The signal-observation pair in \cref{lem:PIctstime} and associated conditions are sufficient but more than likely not necessary to ensure well-posedness of \cref{eq:poiss_phi_mehta}.  For instance, it is well-known that log-concavity is not a strict requirement for the existence of a Poincar\'{e} inequality, as also demonstrated in the next section.  A perturbation style argument used in \cref{lem:PIctsdisctime} poses difficulties for the continuous time observation setting due to the unbounded variation of Brownian motion, unless the signal is of the form $d\hat{X}_t = 0$ (see \cite{Laugesen2015}).

 \begin{remark}
 	\label{rem:vhcond}
 	\cite{VanHandelthesis2007} establishes log-concavity of the posterior measure for the system \cref{eq:sigprocess}-\cref{eq:obsprocess} in the context of filter stability by considering log-concavity of the value function in the relevant stochastic control problem.  We note here that we are able to obtain slightly weaker conditions on the system (cf. Proposition 4.3.6 in \cite{VanHandelthesis2007}), specifically we do not require $c_g \geq 2c_r^2$.	
 \end{remark}

 Well-posedness of the Poisson equation arising in the Crisan \& Xiong filter \cref{eq:pde_beta_crisan} for $d \geq 2$, $\rho_t$ locally bounded and $p$-integrable for some $p \in (1,d)$ and $h \in L^\infty$ was established in \cite{sr:crisan10}.  It is possible to broaden this class of signal and observation drift functions.  It holds from classical regularity theory that when the r.h.s of \cref{eq:pde_beta_crisan} is uniformly bounded (which then implies $(h-\bar{h}) \in L_{\rho_t}^1(\mathbb{R}^{d})$) and $C^2$ that there exists a unique $\beta_t \in C^2$ up to additive constants, so that $\nabla \beta_t$ is unique.  One such example (though there are many) is the system \cref{eq:sigprocess}-\cref{eq:obsprocess} with conditions in \cref{lem:PIctstime}.  
 This is of course conditional on the convergence of \cref{eq:fokkerplankmv}  to the Kushner-Stratonovich equation as $\delta \rightarrow 0$ for the system \cref{eq:sigprocess}-\cref{eq:obsprocess}, which is still an open question.    Analysis of the Poisson equations arising in the Reich filter is difficult without the specification of the mass matrix $M$.

 \section{Well-posedness of the Poisson equation for $\delta > 0$} 
 \label{sec:wellposedctsdisc}
 
 As noted in \cite{Whiteley2018}, an Ornstein-Uhlenbeck signal process combined with a discrete-time observation process that gives rise to log-concave likelihood functions will lead to log-concave posterior densities.  Well-posedness of the FPF Poisson equations \cref{eq:poissphidisc} and \cref{eq:poisspsidisc} then follows from similar arguments as in \cref{theo:wellposedcts}.  
 In this section, we exploit the bounded variation of the $\delta$-approximation to the observation process to establish further conditions that guarantee well-posedness.  In particular, the main result is the following lemma.  The proof appears in \cref{sec:prooflemma51}.   
 
 \begin{lemma}
 	\label{lem:PIctsdisctime}
 	Consider a probability space $(\Sigma, \mathcal{F}, P)$ with the usual conditions.
 	The signal process  $\mathcal{X}_s$ defined on this space evolves according to \cref{eq:sigcts} where $\mathcal{M}: \mathbb{R}^{d} \rightarrow \mathbb{R}^{d}$ is globally Lipschitz with Lipschitz constant $L_\mathcal{M}$.  Observations are generated by the process \cref{eq:obscts} with $h: \mathbb{R}^d \rightarrow \mathbb{R}$ and $h, \nabla h \in L^{\infty}$.  Suppose the observation path is approximated piecewise linearly in time according to \cref{eq:smoothobspath} with mesh size $\delta$.  The initial state $\mathcal{X}_0 \sim \mu_0(x)$ where $\mu_0$ has a density $\rho_0$ with respect to the Lebesgue measure.  
 	Assume that $\rho_0(x)$ satisfies a Poincar\'{e}  inequality with constant $\kappa_0$. Then for all $t \in [0, T]$, with $T = \delta N, N \in \mathbb{Z}_{>0}$, $\rho_{t}^\delta$ which is the solution of \cref{eq:fokkerplankmv} satisfies a Poincar\'{e}  inequality, 
 	\begin{align}
 	\label{eq:PIcontdisc}
 	\int |f(x)|^2 \rho^\delta_{t}(x) dx \leq {\kappa_T} \int |\nabla f(x)|^2 \rho^\delta_{t}(x) dx, 
 	\end{align}    
 	for all test functions $f \in H_{\rho_t^\delta, 0}^1(\mathbb{R}^{d})$ and constant
 	\begin{align*}
 	\kappa_T = (\kappa_0 + T) \exp \left( (2L_\mathcal{M} + |h^2|_\infty)T + 2|h|_\infty \left(\frac{T}{\delta} \right)^2\sup_{n \leq N}|Z_{t_{n+1}} - Z_{t_n}|   \right).
 	\end{align*}
 	
 \end{lemma}
 
  It is clear from the \cref{lem:PIctsdisctime} that unlike in the continuous time setting, the Poincar\'{e} constant grows as the number of observation samples increases.  In this case, we are able to consider a wider class of signal-observation drift functions as compared to the continuous time setting, but at the price of uniformity in time of the well-posedness result.  We then have the following well-posedness result along the same lines as \cref{theo:wellposedcts}. 
 
 \begin{theorem}
 	\label{theo:wellposedctsdisc}
 	\textbf{Well-posedness of FPF Poisson equation for $\delta > 0$.} Consider the signal-observation pair and conditions specified in \cref{lem:PIctsdisctime}.  Then there exists a unique weak solution $\phi \in H^1_{\rho_t^\delta, 0}(\mathbb{R}^{d})$ that satisfies the Poisson equation \cref{eq:poissphidisc} i.e.,  
 	\begin{align}
 	\label{eq:poissweakctsdisc}
 	\int \nabla \phi(x) \cdot \nabla \nu(x) \rho_t^\delta (x) dx &= \int (h(x) - \overline{h}_t) \nu(x) \rho_t^\delta (x)  dx  
 	\rho_t^\delta dx, 
 	\end{align}
 	and there exists a unique weak solution $\psi \in H^1_{\rho_t^\delta, 0}(\mathbb{R}^{d})$ of \cref{eq:poisspsidisc}, i.e. 
 	\begin{align}
 	\label{ew:poissweakctsdisc2}
 	\int \nabla \psi(x) \cdot \nabla \nu(x) \rho_t^\delta (x)  dx &= \int \left( \nabla h^T \nabla \phi(x, \rho_t^\delta) + (\overline{h}_t)^2 - \overline{h_t^2} \right) \rho_t^\delta (x) \nu(x) dx,
 	\end{align}
 	for all test functions $\nu \in H^1_{\rho_t^\delta , 0}(\mathbb{R}^{d})$.
 	
 	\begin{proof}
 		The proof is identical to that of \cref{theo:wellposedcts} except that now 
 		\begin{align*}
 		|T_\phi(\nu)| &:= \left|\int (h(x) - \overline{h}_t)\nu(x) \rho_t^\delta (x) dx \right|,
 		\end{align*}			
 		which is still a bounded linear functional due to \cref{lem:PIctsdisctime} and $h \in L^\infty$.  Additionally, we have 
 		\begin{align*}
 		|T_\psi(\nu)| &:= \left|  \int \left( \nabla h^T \nabla \phi(x, \rho_t^\delta) + (\overline{h}_t)^2 - \overline{h_t^2} \right) \rho_t^\delta (x) \nu(x) dx \right| \\
 		& \leq c \left[|\nabla h|_\infty \left( \int |\nabla \phi|^2 \rho_t^\delta dx  \right)^{1/2} +  2|h|^2_\infty  \right] \left(\int | \nabla \nu|^2 \right)^{1/2},
 		\end{align*}
 		which again is a bounded linear functional due to \cref{lem:PIctsdisctime}, and also $h, \nabla h \in L^\infty$ and $\phi \in H^1_{\rho_t^\delta,0}(\mathbb{R}^{d})$. 
 	\end{proof}
 	
 \end{theorem}
 
 \begin{remark}
 	In both the continuous time and continuous-discrete time setting, the results of \cref{theo:wellposedcts} and \cref{theo:wellposedctsdisc} extend easily to the multivariate observation case by considering the componentwise Poisson equations independently. 
 \end{remark}

 \section{Proofs}
 \label{sec:proofs}
 
 \subsection{Proof of \cref{lem:PIctstime}}
 \label{sec:proofPIlemmacts}
 The basic idea of the proof is to show that when the initial density is log-concave, the posterior density stays log-concave under the signal dynamics and likelihood transformation.  We start by considering a time discretisation of the Kallianpur-Striebel formula, which is useful in this context because of the Girsanov transformation which then allows us to work with integrals against a Wiener measure.  By the process of induction and Brascamp-Lieb type inequalities, we then show that the time discretised posterior is log-concave, and that this property is maintained in the continuous time limit. \\
 \\
 Recall that in Lemma \ref{lem:kushnerrho} we showed that the conditional probability density of $X_t$ given $\mathcal{Z}_t$ evolves according to the Kushner-Stratonovich equation, which implies that $\rho_t = \theta_t$.  We can therefore work with the well-known equations in filtering theory to establish the required Poincar\'{e} inequality.  Specifically, the continuous time posterior expectation at some arbitrary time $T$ is given by the Kallianpur-Striebel formula 
 \begin{align}
 	\label{eq:kallinapurstriebel}
 	\mathbb{E}_P(f(\mathcal{X}_T) | \mathcal{Z}_T)(\omega) = \frac{\mathbb{E}_{\tilde{P}}(f(\mathcal{X}_T) M_T(\mathcal{X}, Z(\omega))}{\mathbb{E}_{\tilde{P}}( M_T(\mathcal{X}, Z(\omega)))}, 
 \end{align}
 with 
 \begin{align*}
 	M_T(\mathcal{X}, Z) & = \exp \left(\int_0^T H\mathcal{X}_s dZ_s - \frac{1}{2} \int_0^T |H \mathcal{X}_s|^2  ds  \right) = \frac{dP}{d\tilde{P}}, 	
 \end{align*}
 where $\tilde{\mathbb{E}}$ indicates the expectation with respect to the probability measure $\tilde{P}$, $Z$ is a Wiener process under $\tilde{P}$ and $\mathcal{X}$ has the same law under $P$ and $\tilde{P}$. The Kallianpur-Striebel formula can also be reformulated as 
 \begin{align}
 	\label{eq:KallStriebpath}
 	\mathbb{E}_P(f(\mathcal{X}_T) | \mathcal{Z}_T) (\omega) = \frac{\int_{C[0,T]} f(x_T) M_T(x, Z(\omega)) \mu_{\mathcal{X}}(dx)  }{\int_{C[0,T]}  M_T(x, Z(\omega)) \mu_{\mathcal{X}}(dx) },
 \end{align}
 where $\mu_{\mathcal{X}}(dx)$ indicates the path measure induced by the signal process $\mathcal{X}$ and $C_{[0,T]}$ is the space of continuous functions on the time interval $[0,T]$.  Applying Girsanov and using the fact that the law of $\mathcal{X}$ is unchanged under $P$ or $\tilde{P}$, we have 
 \begin{align}
 	\label{eq:expecwiener}
 	\mathbb{E}_{\tilde{P}} (g(\mathcal{X}, Z(\omega))) = \int_{C[0,T]}  f(v_T) N_T(v,Z(\omega)) \mu_V(dv),
 \end{align}
 where 
 \begin{align*}
 	g(\mathcal{X}, Z(\omega)) & = f(\mathcal{X}_T) M_T(\mathcal{X}, Z(\omega)), \\
 	N_T(V, Z(\omega)) & = M_T(V,Z(\omega)) \exp \left( \int_0^T \nabla U (V_s) \cdot dV_s - \frac{1}{2} \int_0^T |\nabla U(V_s)|^2 ds  \right),
 \end{align*}
 and $\mu_V(dv)$ indicates the Wiener measure on path space.  In order to apply Girsanov, we must confirm the Novikov condition, i.e. , 
 \begin{align*}
 	\mathbb{E}_{\tilde{P}} \left( \exp \left( \frac{1}{2} \int_0^T |\nabla U (V_s)|^2 ds  \right)  \right) < \infty,
 \end{align*}
 which holds due to \ref{cond:a4}. It\^{o}'s formula gives 
 \begin{align*}
 	\int_0^T \nabla U (V_s) \cdot dV_s = U(V_T) - U(V_0) - \frac{1}{2} \int_0^T \Delta U (V_s) ds, 
 \end{align*}
 which then leads to 
 \begin{align*}
 	&N_T(V, Z(\omega)) = \\
 	&\exp \left( U(V_T) - U(V_0) + \int_0^T HV_s dZ_s(\omega) - \frac{1}{2} \int_0^T |H V_s|^2 + (\Delta U + |\nabla U|^2)(V_s) ds  \right). 
 \end{align*}
 
 Now consider the following time discretisation of the conditional expectation.  Specifically, a finite time interval $[0,T]$ is discretised into the following sequence $\left \lbrace t_0, t_1, ..., t_N \right \rbrace$ of sampling time instances with $t_0 = 0 < t_1 < ... < t_N = T$, and time increment $\Delta t = t_i - t_{i-1}$ for all $i = 1, ..., N$, where $\Delta t = \frac{T}{N}$. 
 Define 
 \begin{align*}
 	\mathcal{N}_{N} := \prod_{i=0}^{N-1} \exp \left( HV_{t_i} \Delta Z_{i+1} - \frac{1}{2}  \left[|HV_{t_i}|^2 + (\Delta U + |\nabla U|^2)(V_{t_i})  \right]\Delta t \right),
 \end{align*} 
 where $\Delta Z_{{i+1}}:= Z_{t_{i+1}} - Z_{t_i}$, and we have dropped the $\omega$ for notational ease.  We can now work with the time-discretised version to establish a Poincar\'{e} inequality in both continuous and discrete time, for the system \cref{eq:sigprocess}- \cref{eq:obsprocess}.  That is, by using the finite dimensional distributions of Brownian motion we have that \cref{eq:expecwiener} can be expressed as 
 \begin{align*}
 	& {\mathbb{E}}_{\tilde{P}} (f(v_T) \exp(U(v_T) - U(v_0)) \mathcal{N}_{N}) \\
 	& = \int_{C[0,T]}  (f(v_T) \exp(U(v_T) - U(v_0)) \prod_{i=0}^{N-1} \exp \left( Hv_{t_i} \Delta Z_{i+1} - \frac{1}{2}  R(v_{t_i})\Delta t \right) \mu^V(dv) \\
 	& = \int f(v_T)  \exp(U(v_T)) \int Q(v_{t_{N-1}}) q_{\Delta t}(v_{t_{N-1}}, v_{t_{N}})  \int Q(v_{t_{N-2}}) q_{\Delta t}(v_{t_{N-2}}, v_{t_{N-1}}) \cdots \\
 	& \int Q(v_{{0}}) q_{\Delta t}(v_{{0}}, v_{t_{1}}) \exp(-U(v_0)) \rho_0(v_0) dv_0 dv_{t_{1}} \cdots dv_{t_{N-2}} dv_{t_{N-1}},
 \end{align*}
 where 
 \begin{align}
 	\nonumber
 	Q(v_{t_{i}}) &:= \exp \left( Hv_{t_i} \Delta Z_{i+1} - \frac{1}{2} R(v_{t_i}) \Delta t \right), \\
 	\label{eq:funcR}
 	R(v) &:=  |Hv|^2 + (\Delta U + |\nabla U|^2)(v) , \\
 	\nonumber
 	q_{\Delta t}(v_{t_{i-1}}, v_{t_i}) &:= \frac{1}{\sqrt{2 \pi \Delta t}} \exp\left(-\frac{|v_{t_i} - v_{t_{i-1}}|^2}{2 \Delta t} \right),
 \end{align}
 from which we obtain the following recursion: 
 \begin{align}
 	\nonumber
 	\check{\rho}_0(v_0) &=  \exp(-U(v_0)) \rho_0(v_0):= \exp(-\check{\mathcal{G}}_0(v_0)), \\
 	\label{eq:rhorecursion}
 	\check{\rho}_i(v_{t_i}) &=  \int Q(v_{t_{i-1}}) q_{\Delta t}(v_{t_{i-1}}, v_{t_i}) \check{\rho}_{i-1}(v_{t_{i-1}}) dv_{t_{i-1}} \quad \forall \enskip i = 1, 2, \cdots, N.
 \end{align}
 By \ref{cond:a2}, we have that $\check{\mathcal{G}}_0(v_0)$ is strongly convex with parameter $c_g$.  Then consider 
 \begin{align*}
 	\check{\rho}_1(v_{t_1}) &\propto  \int \exp \left( Hv_{0} \Delta Z_{1} - \frac{1}{2} R(v_{0}) \Delta t -\frac{|v_{t_1} - v_{0}|^2}{2 \Delta t} - \check{\mathcal{G}}_0(v_0) \right)   dv_{0} \\
 	& =: \int \exp(f(v_0, v_{t_1})) dv_0 \\
 	& =: \exp(-\check{\mathcal{G}}_1(v_{t_1})),
 \end{align*}
 and by Theorem \ref{theo:prekopa}, we have that 
 \begin{align*}
 	\nabla_{v_{t_1}}^2 \check{\mathcal{G}}_1 & \succcurlyeq \frac{\int (\nabla_{v_{t_1}}^2 f - \nabla_{v_0 v_{t_1}}^2 f (\nabla_{v_0}^2 f)^{-1} \nabla_{v_{t_1} v_0}^2 f  ) e^{-f(v_0,v_{t_1})} dv_0}{\int e^{-f(v_0,v_{t_1})}dv_0 } \\
 	&= \int \left(\frac{1}{ \Delta t}I - \frac{1}{( \Delta t)^2}  \left( \frac{\Delta t}{2} \nabla_{v_0}^2 R + \nabla_{v_0}^2 g_0 + \frac{1}{\Delta t} I \right)^{-1}   \right) e^{-f(v_0,v_{t_1})} dv_0 \Big/ \int e^{-f(v_0,v_{t_1})}dv_0 \\
 	& \succcurlyeq \left( \frac{1}{\Delta t} - \frac{1}{(\Delta t)^2}  \left( c_g + \frac{\Delta t c_r}{2} + \frac{1}{\Delta t} \right)^{-1} \right) I \\
 	& =  \left( \frac{c_g + \frac{\Delta t c_r}{2}}{1+ \Delta t \left(c_g + \frac{\Delta t c_r}{2} \right)} \right) I \\
 	& =: \gamma_1 I.
 \end{align*}
 By induction, one obtains that $\check{\rho_i} (v_{t_i}) \propto \exp(-\check{\mathcal{G}}_i(v_{t_i}))$ with $\nabla_{v_{t_i}}^2 \check{\mathcal{G}}_i \succcurlyeq \gamma_i I $ (i.e.,~it is log-concave) for all $ i = 1, 2, \cdots N$, where 
 \begin{align*}
 	\gamma_i  = m(\gamma_{i-1})  &=  \left( \frac{\gamma_{i-1} + \frac{\Delta t c_r}{2}}{1+ \Delta t \left(\gamma_{i-1} + \frac{\Delta t c_r}{2} \right)} \right), \\ 
 	\gamma_0 & = c_g. 
 \end{align*} 
 Since $0 < \frac{d m(x)}{dx} < 1$ for $c_g, c_r > 0$ and $\Delta t > 0$, it follows that $\gamma_i$ approaches a stable fixed point $\gamma_*$ monotonically in time, where 
 \begin{align*}
 	\gamma_* &= \frac{- \Delta t c_r + \sqrt{(\Delta t {c_r})^2 + 8c_r}}{4}, 
 \end{align*}
 and 
 \begin{align*}
 	\gamma_* \rightarrow \sqrt{\frac{ c_r}{2}} \quad \text{as} \enskip  \Delta t \rightarrow 0,
 \end{align*}
 which together implies that 
 \begin{align*}
 	\gamma_i \geq \min \left( c_g, \sqrt{\frac{c_r}{2}}  \right) \quad \forall \enskip i = 0, 1, \cdots, N.  
 \end{align*}
That is, we have shown that $\check{\rho}_i(v_{t_i})$ is log-concave for all  $i=0, 1, \cdots, N$.   It is clear that the time discretised posterior density at time $T$ is given by 
\begin{align}
	\label{eq:defndiscpost}
	{\rho}_N(v_{t_N})  \propto  \exp(U(v_{t_N}))   \check{\rho}_N (v_{t_N}),
\end{align}  
 where we use the integer subscript $N$ to differentiate from the continuous time posterior $\rho_T$.  Furthermore, due to \ref{cond:a1}, ${\rho}_N(v_{t_N})$ is log-concave with parameter $c_u + \min \left( c_g, \sqrt{\frac{c_r}{2}}  \right)$, which holds even as $N \rightarrow \infty$.  By a direct application of the Brascamp-Lieb inequality (Theorem \ref{theo:brascamplieb}), we have that  
 \begin{align*}
 	\int |f(v_{t_N})|^2 {\rho}_N(v_{t_N}) dv_{t_N} \leq \frac{1}{c_\gamma} \int |\nabla f(v_{t_N})|^2 {\rho}_N(v_{t_N}) dv_{t_N} 
 \end{align*}    
 with $c_\gamma := c_u + \min \left( \omega, \sqrt{\frac{c_r}{2}}  \right)$ for all test functions $f \in H^1_{{\rho}_N,0}(\mathbb{R}^{d})$,  which can be equivalently stated as 
 \begin{align}
 	\label{eq:PIdisclogconc}
 	\mathbb{E}_{\tilde{P}} (|f(v_T)|^2\exp(U(v_T) - U(v_0)) \mathcal{N}_N ) \leq \frac{1}{c_\gamma} \mathbb{E}_{\tilde{P}} (|\nabla f(v_T)|^2\exp(U(v_T) - U(v_0)) \mathcal{N}_N )
 \end{align} 
 for all $N \in \mathbb{N}_{>0}$.  Finally, from the definition of an It\^{o}  stochastic integral we have that 
 \begin{align*}
 	\mathcal{N}_N \xrightarrow{L^2(\tilde{P})} N_T \quad \text{as} \quad N \rightarrow \infty
 \end{align*}
 and since $L^2$ convergence implies weak convergence, it follows that 
 \begin{align*}
 	\lim_{N \rightarrow \infty} {\mathbb{E}}_{\tilde{P}} (g(V_T) \exp(U(V_T) - U(V_0)) \mathcal{N}_{N}) = \mathbb{E}_{\tilde{P}}(g(V_T) \exp(U(V_T) - U(V_0)) N_T),
 \end{align*}
 which together with \cref{eq:PIdisclogconc} implies \cref{eq:PIlogconcont}. 
 
 \qed
 
 \subsection{Proof of  \cref{lem:PIctsdisctime}}	
 \label{sec:prooflemma51}
 We proceed by first establishing a Poincar\'{e}  inequality for the law of the signal process, making use of the fact that 
 \begin{align}
 \label{eq:marginalexp}
 \int_{C_{[0,T]}} |f(x_T)|^2 \mu_{\mathcal{X}}(dx) = \int_{\mathbb{R}^d} |f(x_T)|^2 \mu_{\mathcal{X}_T}(dx_T),
 \end{align}
 where $\mu_{\mathcal{X}_T}$ is the probability measure induced by the solution of \cref{eq:sigcts} at time $T$ under $P$.  Under the conditions on $\mathcal{M}$ and $h$, $\mu_{\mathcal{X}_T}$ is absolutely continuous with respect to the Lebesgue measure, and we denote its density by $\rho_T^-$.  Consider the following Euler-Maruyama discretisation of the signal process over $[0, T]$ with time step $\Delta \tau = \frac{T}{L}$: 
 \begin{equation}
 X_l = X_{l-1} + \mathcal{M}(X_{l-1}) \Delta \tau +  \sqrt{\Delta \tau}\Delta V_l,
 \end{equation}
 where $\Delta V_l \sim N(0,I)$.  Denote by $\check{\rho}_l$ the probability density of $X_{l} \enskip \forall \enskip l \in \{1, 2, \cdots, L\}$.  The following lemma will be of use, which appears in the literature in various forms and we present here the precise form needed for our purposes.  
 
 \begin{lemma}
 	\label{lem:invPIlip}
 	\textbf{Invariance of Poincar\'{e} inequality under Lipschitz Transformations}
 	Let $\eta$ and $\nu$ be two density functions associated to the random variables $X$ and $Y$, respectively, and let  $\Theta: \mathbb{R}^d \rightarrow \mathbb{R}^d$ be a Lipschitz  map with Lipschitz constant $L_\Theta$ such that $Y = \Theta(X)$.  Suppose $\eta$ satisfies a Poincar\'{e} inequality 
 	\begin{equation}
 	\label{eq:pimu}
 	\exists \enskip c > 0 \enskip | \enskip \forall \enskip f \in H^1_{\eta, 0}(\mathbb{R}^d), \enskip  \enskip \int |f(x)|^2 \eta(x) dx \leq c \int |\nabla f(x)|^2 \eta(x) dx.
 	\end{equation}
 	Then $\nu$ satisfies a PI also:
 	\begin{equation*}
 	\exists \enskip c_2 > 0 \enskip | \enskip \forall \enskip f \in H^1_{\nu,0}(\mathbb{R}^d), \enskip  \enskip \int |f(y)|^2 \nu(y) dy \leq c_2 \int |\nabla f(y)|^2 \nu(y) dy,
 	\end{equation*}
 	with $c_2 = c L_\Theta.$ 
 	
 	\begin{proof}
 		By change of variables,
 		\begin{equation*}
 		\int \psi(y) \nu(y) dy = \int \psi(\Theta(x)) \eta(x) dx,
 		\end{equation*}
 		for all integrable $\psi$.
 		Let $\psi(y) = |f(y)|^2$, then using \cref{eq:pimu} we have 
 		\begin{subequations}
 			\begin{align*}
 			\int |f(y)|^2 \nu(y) dy &= \int |f(\Theta(x))|^2 \eta(x) dx \\
 			&\leq c\int |\nabla [f(\Theta(x))]|^2 \eta(x) dx  \\
 			& \leq c\int |\nabla [f(\Theta)]|^2  |\nabla \Theta(x)|^2 \eta(x) dx \\
 			& \leq c L_\Theta  \int |\nabla [f(\Theta)]|^2 \eta(x) dx \\
 			& = c L_\Theta \int |\nabla f(y)|^2 \nu(y) dy,
 			\end{align*}
 		\end{subequations}	
 		where in the above $|.|$ refers to the $L^2$ induced matrix norm.
 	\end{proof}
 \end{lemma}
 
 Since $\mathcal{M}$ is globally Lipschitz, by a direct application of \cref{lem:invPIlip} and the fact that Poincar\'{e} inequalities are stable under convolutions (Corollary 3.1 in \cite{Chafai2004}) we have that $\check{\rho}_{1}$ satisfies a Poincar\'{e} inequality with constant 
 \begin{align*}
 \check{\kappa}_{1} &=  \kappa_{0}(1+L_\mathcal{M} \Delta \tau)^2 + \Delta \tau. 
 \end{align*} 
 By induction, the Poincar\'{e} constant for $\check{\rho}_{ L } (x)$ is given by 
 \begin{align*}
 \check{\kappa}_{ L} = \alpha^L \kappa_0 + \beta \sum_{j=0}^{L-1} \alpha^j,
 \end{align*}
 where $\alpha = (1 + L_\mathcal{M} \Delta \tau)^2$ and $\beta = \Delta \tau$.  This implies that 
 \begin{align}
 \label{eq:PIsigcontdisc}
 \int |f(x)|^2 \check{\rho}_{ L } (x) dx \leq \check{\kappa}_L \int |\nabla f(x)|^2 \check{\rho}_{ L } (x) dx
 \end{align}
 for $f \in H^1_{\check{\rho}_L, 0}(\mathbb{R}^{d})$.  Now, it is well-known that the Euler-Maruyama scheme converges in $L^1$, which therefore also implies weak convergence, i.e., 
 \begin{align}
 \label{eq:eulerweakconv}
 \lim_{L \rightarrow \infty} \mathbb{E}_P(X_L) = \mathbb{E}_P(\mathcal{X}_{t}).
 \end{align}
 Furthermore, we have that 
 \begin{align*}
 \lim_{L \rightarrow \infty}  {\check{\kappa}_L} \int |\nabla f|^2 \check{\rho}_{ L } (x) dx = \left(  \lim_{L \rightarrow \infty}  {\check{\kappa}_L} \right)  \left( \lim_{L \rightarrow \infty}  \int |\nabla f|^2 \check{\rho}_{ L } (x) dx  \right),
 \end{align*}
 since both limits exist, and the first limit on the r.h.s is given by ${\kappa_T^-}$ where 
 \begin{align*}
 {\kappa}^-_T =& \lim_{L \rightarrow \infty} \check{\kappa}_L \\
 =&  \lim_{L \rightarrow \infty}  \left[ \kappa_0 \left(1+ \frac{L_\mathcal{M} T}{L}\right)^{2l}  + \frac{T}{L} \sum_{j=0}^{L-1} \left(1+ \frac{L_\mathcal{M} T}{L}\right)^{2j}  \right ].
 \end{align*}
 It follows that 
 \begin{equation*}
 \lim_{L \rightarrow \infty}  \kappa_0 \left(1+ \frac{L_\mathcal{M} T}{L}\right)^{2l} = \kappa_0 e^{2L_\mathcal{M} T}
 \end{equation*}
 and 
 \begin{equation*}
 \sum_{j=0}^{L-1} \left(1+ \frac{L_\mathcal{M} T}{L}\right)^{2j} \leq Le^{2L_\mathcal{M} T},
 \end{equation*}
 so that
 \begin{equation*}
 \lim_{L \rightarrow \infty} \frac{T}{L} \sum_{j=0}^{L-1} \left(1+ \frac{L_\mathcal{M} T}{L}\right)^{2j} \leq \lim_{L \rightarrow \infty} \frac{T}{L} L e^{2L_\mathcal{M} T} = Te^{2L_\mathcal{M} T}.
 \end{equation*}
 Combining the above, we have
 \begin{equation}
 \label{eq:dissignalPC}
 \kappa^-_T = \exp({2L_\mathcal{M} T})(\kappa_0 + T).
 \end{equation}
 Finally, \cref{eq:dissignalPC}, \cref{eq:eulerweakconv} and \cref{eq:marginalexp} together imply that  
 \begin{equation}
 \label{eq:PIsignal}
 \int_{C_{[0,T]}} |f(x_T)|^2 \mu_{\mathcal{X}}(dx) \leq {\kappa_T^-} \int_{C_{[0,T]}} |\nabla f(x_T)|^2 \mu_{\mathcal{X}}(dx).
 \end{equation}	
 As previously, since we have $\rho_t^\delta = \theta_t^\delta$ by construction, we can make use of the $\delta$-approximation of the Kallianpur-Striebel formula \cref{eq:KallStriebpath}.  We note that the smooth approximation to the exponential martingale, 
 \begin{align*}
 M_T(\mathcal{X},Z^\delta(\omega)) & = \exp \left( \int_0^T h(\mathcal{X}_s) dZ^\delta_s(\omega)- \frac{1}{2} \int_0^T |h(\mathcal{X}_s)|^2 ds  \right), \\
 & = \exp \left( \int_0^T h(\mathcal{X}_s) \frac{Z_{t_{n+1}}(\omega) - Z_{t_n}(\omega)}{\delta} ds - \frac{1}{2} \int_0^T |h(\mathcal{X}_s)|^2 ds  \right).
 \end{align*}
 is uniformly bounded when $h$ is bounded, since $Z_s^\delta$ is of bounded variation.  Then by the classical perturbation result due to Holley \& Stroock we have that   
 \begin{align*}
 \int_{C_{[0,T]}} |f(x_T)|^2 M_T(x, Z^\delta(\omega)) \mu_{\mathcal{X}}(dx) \leq {\kappa_T} \int_{C_{[0,T]}} |\nabla f(x_T)|^2 M_T(x, Z^\delta(\omega))\mu_{\mathcal{X}}(dx).
 \end{align*}
 where $\kappa_T = \kappa_T^- \exp(\text{osc}(M_T(\mathcal{X}, Z^\delta)))$.  Furthermore, we have that 
 \begin{align*}
 \exp(\text{osc}(M_T(\mathcal{X}, Z^\delta(\omega)))) \leq \exp \left( T|h^2|_\infty + 2|h|_\infty\left( \frac{T}{\delta} \right)^2\sup_{n \leq N}|Z_{t_{n+1}}(\omega) - Z_{t_n}(\omega)|  \right),
 \end{align*}
 which gives the desired result. 
 
 \qed

 \section{Conclusions}	
 %
 
 A number of particle filters have been proposed over the last couple of decades with the common feature that the update step is governed by a type of control law.  This feature makes them an attractive alternative to traditional sequential Monte Carlo which scales poorly with the state dimension due to weight degeneracy.  This article proposes a unifying framework that allows to systematically derive the McKean-Vlasov representations of these filters for the discrete time and continuous time observation case, taking inspiration from the smooth approximation of the data in \cite{sr:crisan10} and \cite{Clark2005}.  We highlight the various choices on the coefficients of the process that then leads to the different filters considered here, i.e.,~the FPF \cite{Yang2011}, \cite{Yang2014}, \cite{Laugesen2015}, the so-called Crisan \& Xiong filter \cite{sr:crisan10} and the so-called Reich filter \cite{sr:reich10}.  The limiting form of the Crisan \& Xiong and Reich filters is derived, thereby extending them to the continuous time observation setting. Additionally, this framework leads to formulations of the FPF and Reich filters for the continuous signal-discrete observation setting that avoids the need for a two step predict-update procedure (the so-called $\delta$-FPF and $\delta$-Reich filters, respectively).  All filters require the solution of a Poisson equation defined on $\mathbb{R}^{d}$, for which existence and uniqueness of solutions can be a nontrivial issue.  We therefore also establish conditions on the signal-observation system that ensures well-posedness of the weighted Poisson equation arising in the FPF.      
 
 There are several future research directions arising from this work that could be investigated.  For instance, there are many other filters (e.g.,~\cite{daumhuangparticleflow}) that could be further analysed under this framework as well, in particular to obtain their limiting forms.  Most importantly, we have only provided representations of the filters for the continuous time observation limit, and a rigorous justification of the existence of such processes as well as convergence will be investigated in a future publication.  There also exist possibilities to consider relaxations on the conditions on $\mathcal{M}$ and $h$ to ensure well-posedness of the Poisson equations arising in these filters.

\bibliographystyle{siamplain}
\bibliography{literature_PIpaper}    

\newpage 

\appendix
\section{Some useful theorems}
\begin{theorem}
\label{theo:brascamplieb}	
\textbf{Brascamp-Lieb inequality \cite{Brascamp1976}}. Let $\rho(x) = \exp(-\mathcal{G}(x))$ where $x \in \mathbb{R}^d$ and $\mathcal{G} \in C^2(\mathbb{R}^d)$ are strictly convex. Then it holds for all $f(x) \in H^2_\rho(\mathbb{R}^d)$ that 
\begin{align*}
	\int |f(x)|^2 \rho(x)dx - \left( \int f(x) \rho(x) dx  \right)^2 \leq \int \langle \nabla f, [\nabla^2 \mathcal{G}]^{-1} \nabla f \rangle \rho(x)dx.
\end{align*}	    
\end{theorem}	

As detailed in \cite{Brascamp1976}, a useful corollary of the above theorem is the following (see Theorem 4.2 in \cite{Brascamp1976})
\begin{theorem}
	\label{theo:prekopa}
	Let $\mu(x) = \mu(y,z) = \exp(-\mathcal{G}(y,z))$ with $y \in \mathbb{R}^m, z \in \mathbb{R}^n$ and $d = m+n$.  Suppose $\mathcal{G}$ satisfies the assumptions of Theorem \ref{theo:brascamplieb}. Consider $\rho(y) = \exp(-\mathcal{H}(y)) = \int \exp(-\mathcal{G}(y,z)) dz$.  Then it holds that  
	\begin{align*}
		\nabla_{y}^2 \mathcal{H}  \succcurlyeq \frac{\int (\nabla_{y}^2 \mathcal{G} - \nabla_{yz}^2 \mathcal{G} (\nabla_{z}^2 \mathcal{G})^{-1} \nabla_{zy}^2 \mathcal{G} ) \mu(y,z) dz}{\int \mu(y,z) dz }. 
	\end{align*}	
\end{theorem}

\end{document}

%% file: ex_shared.tex
\usepackage{lipsum}
\usepackage{amsfonts}
\usepackage{amssymb}
\usepackage{graphicx}
\usepackage{epstopdf}
\usepackage{algorithmic}
 \usepackage{array}
\newcolumntype{M}[1]{>{\centering\arraybackslash}m{#1}}
\newcolumntype{N}{@{}m{0pt}@{}}
\usepackage{enumitem} 
\ifpdf
  \DeclareGraphicsExtensions{.eps,.pdf,.png,.jpg}
\else
  \DeclareGraphicsExtensions{.eps}
\fi

\newsiamremark{remark}{Remark}
\newsiamremark{assumption}{Assumption}
\newsiamremark{hypothesis}{Hypothesis}
\crefname{hypothesis}{Hypothesis}{Hypotheses}
\newsiamthm{claim}{Claim}

\newcommand{\norm}[1]{\left\lvert#1\right\rvert}

\headers{McKean-Vlasov SDEs in Nonlinear Filtering}{S. Pathiraja, S. Reich, and W. Stannat}

\title{McKean-Vlasov SDEs in nonlinear filtering\thanks{Submitted to the editors DATE.
\funding{This research has been partially funded by Deutsche Forschungsgemeinschaft (DFG) - Project-ID 318763901 - SFB1294.}}}

\author{Sahani Pathiraja\thanks{Universit\"at Potsdam, 
		Institut f\"ur Mathematik, Karl-Liebknecht-Str. 24/25, D-14476 Potsdam, Germany 
  (\email{pathiraja@uni-potsdam.de}).}
\and Sebastian Reich\thanks{Universit\"at Potsdam, Institut f\"ur Mathematik, Karl-Liebknecht-Str. 24/25, D-14476 Potsdam, Germany
  (\email{sreich@math.uni-potsdam.de}).}
\and Wilhelm Stannat \thanks{TU Berlin, Institut f\"ur Mathematik, Str. des 17. Juni 136, D-10623 Berlin, Germany (\email{stannat@math.tu-berlin.de}) .} }

\usepackage{amsopn}


%% file: SICON_Pathirajaetal_FINAL.bbl
\begin{thebibliography}{10}

\bibitem{sr:crisan}
{\sc A.~Bain and D.~Crisan}, {\em Fundamentals of stochastic filtering},
  vol.~60 of Stochastic modelling and applied probability, Springer-Verlag,
  New-York, 2008.

\bibitem{Bakrygentled}
{\sc D.~Bakry, I.~Gentil, and M.~Ledoux}, {\em {Analysis and Geometry of Markov
  Diffusion Operators}}, Springer, 2014.

\bibitem{Bergemann2012}
{\sc K.~Bergemann and S.~Reich}, {\em {An ensemble Kalman-Bucy filter for
  continuous data assimilation}}, Meteorologische Zeitschrift, 21 (2012),
  pp.~213--219, \url{https://doi.org/10.1127/0941-2948/2012/0307}.

\bibitem{Bishop2018}
{\sc A.~N. Bishop, P.~{Del Moral}, and S.~D. Pathiraja}, {\em {Perturbations
  and projections of Kalman-Bucy semigroups}}, Stochastic Processes and their
  Applications, 128 (2018), \url{https://doi.org/10.1016/j.spa.2017.10.006}.

\bibitem{Brascamp1976}
{\sc H.~J. Brascamp and E.~H. Lieb}, {\em {On extensions of the Brunn-Minkowski
  and Pr{\'{e}}kopa-Leindler theorems, including inequalities for log concave
  functions, and with an application to the diffusion equation}}, Journal of
  Functional Analysis, 22 (1976), pp.~366--389,
  \url{https://doi.org/10.1016/0022-1236(76)90004-5}.

\bibitem{Chafai2004}
{\sc D.~Chafai}, {\em {Entropies, convexity, and functional inequalities}},
  Journal of Mathematics of Kyoto University, 44 (2004), pp.~325--363,
  \url{https://doi.org/10.1215/kjm/1250283556}.

\bibitem{Clark2005}
{\sc J.~M. Clark and D.~Crisan}, {\em {On a robust version of the integral
  representation formula of nonlinear filtering}}, Probability Theory and
  Related Fields, 133 (2005), pp.~43--56,
  \url{https://doi.org/10.1007/s00440-004-0412-5}.

\bibitem{Coghi2021}
{\sc M.~Coghi, T.~Nilssen, and N.~N{\"{u}}sken}, {\em {Rough McKean-Vlasov
  dynamics for robust ensemble Kalman filtering}},  (2021), pp.~1--41,
  \url{https://arxiv.org/abs/2107.06621}.

\bibitem{sr:crisan10}
{\sc D.~Crisan and J.~Xiong}, {\em {Approximate McKean-Vlasov representations
  for a class of SPDEs}}, Stochastics, 82 (2010),
  \url{https://doi.org/10.1080/17442500902723575}.

\bibitem{daumhunghomot}
{\sc F.~Daum and J.~Huang}, {\em {Nonlinear filters with particle flow induced
  by log-homotopy}}, in Signal Processing, Sensor Fusion, and Target
  Recognition XVIII, I.~Kadar, ed., vol.~7336, International Society for Optics
  and Photonics, SPIE, 2009, pp.~76--87,
  \url{https://doi.org/10.1117/12.814241}.

\bibitem{daumhuangparticleflow}
{\sc F.~Daum, J.~Huang, and A.~Noushin}, {\em {Exact particle flow for
  nonlinear filters}}, in Signal Processing, Sensor Fusion, and Target
  Recognition XIX, I.~Kadar, ed., vol.~7697, International Society for Optics
  and Photonics, SPIE, 2010, pp.~92--110,
  \url{https://doi.org/10.1117/12.839590}.

\bibitem{daumstochpart}
{\sc F.~Daum, J.~Huang, and A.~Noushin}, {\em {New theory and numerical results
  for Gromov's method for stochastic particle flow filters}}, in 21st
  International Conference on Information Fusion (FUSION), 2018, pp.~108--115,
  \url{https://doi.org/10.23919/ICIF.2018.8455287}.

\bibitem{deWiljes2020}
{\sc J.~de~Wiljes and X.~T. Tong}, {\em {Analysis of a localised nonlinear
  ensemble Kalman-Bucy filter with complete and accurate observations}},
  Nonlinearity, 33 (2020), pp.~4752--4782,
  \url{https://doi.org/10.1088/1361-6544/ab8d14}.

\bibitem{DelMoral2021}
{\sc P.~del Moral and E.~Horton}, {\em {A theoretical analysis of
  one-dimensional discrete generation ensemble Kalman particle filters}},
  (2021), \url{https://arxiv.org/abs/2107.01855}.

\bibitem{sr:evensen03}
{\sc G.~Evensen}, {\em {The ensemble Kalman filter: theoretical formulation and
  practical implementation}}, Ocean Dynamics, 53 (2003), pp.~343--367,
  \url{https://doi.org/10.1007/s10236-003-0036-9}.

\bibitem{Hu2002}
{\sc Y.~Hu, G.~Kallianpur, and J.~Xiong}, {\em {An approximation for the Zakai
  equation}}, Applied Mathematics and Optimization, 45 (2002), pp.~23--44,
  \url{https://doi.org/10.1007/s00245-001-0024-8}.

\bibitem{8618878}
{\sc J.-W. Kim, A.~Taghvaei, and P.~G. Mehta}, {\em Derivation and extensions
  of the linear feedback particle filter based on duality formalisms}, in 2018
  IEEE Conference on Decision and Control (CDC), 2018, pp.~7188--7193,
  \url{https://doi.org/10.1109/CDC.2018.8618878}.

\bibitem{Laugesen2015}
{\sc R.~S. Laugesen, P.~G. Mehta, S.~P. Meyn, and M.~Raginsky}, {\em {Poissons
  Equation in Nonlinear Filtering}}, SIAM Journal on Control and Optimization,
  53 (2015), pp.~501--525, \url{https://doi.org/10.1137/13094743X}.

\bibitem{Mitter2003}
{\sc S.~K. Mitter and N.~J. Newton}, {\em {A variational approach to nonlinear
  estimation}}, SIAM Journal on Control and Optimization, 42 (2003),
  pp.~1813--1833, \url{https://doi.org/10.1137/S0363012901393894}.

\bibitem{Nott2011}
{\sc D.~J. Nott, L.~Marshall, and T.~M. Ngoc}, {\em {The ensemble Kalman filter
  is an ABC algorithm}}, Statistics and Computing, 22 (2011), pp.~1273--1276,
  \url{https://doi.org/10.1007/s11222-011-9300-x}.

\bibitem{Nusken2019}
{\sc N.~N{\"{u}}sken, S.~Reich, and P.~J. Rozdeba}, {\em {State and parameter
  estimation from observed signal increments}}, Entropy, 21 (2019), pp.~1--25,
  \url{https://doi.org/10.3390/e21050505}.

\bibitem{Pathiraja2020}
{\sc S.~Pathiraja}, {\em {L2 convergence of smooth approximations of stochastic
  differential equations with unbounded coefficients}},  (2020), pp.~1--15,
  \url{https://arxiv.org/abs/2011.13009}.

\bibitem{sr:reich10}
{\sc S.~Reich}, {\em {A dynamical systems framework for intermittent data
  assimilation}}, BIT Numerical Mathematics, 51 (2010), pp.~235--249,
  \url{https://doi.org/10.1007/s10543-010-0302-4}.

\bibitem{Reich2019}
{\sc S.~Reich}, {\em {Data assimilation: The Schr{\"{o}}dinger perspective}},
  Acta Numerica, 28 (2019), pp.~635--711,
  \url{https://doi.org/10.1017/S0962492919000011}.

\bibitem{Stannat2005}
{\sc W.~Stannat}, {\em {Stability of the filter equation for a time-dependent
  signal on $\mathbb{R}^d$}}, Applied Mathematics and Optimization, 52 (2005),
  pp.~39--71, \url{https://doi.org/10.1007/s00245-005-0820-7}.

\bibitem{Surace2019HowTA}
{\sc S.~C. Surace, A.~Kutschireiter, and J.-P. Pfister}, {\em How to avoid the
  curse of dimensionality: Scalability of particle filters with and without
  importance weights}, SIAM Review, 61 (2019), pp.~79--91,
  \url{https://doi.org/10.1137/17M1125340}.

\bibitem{Taghvaei2017}
{\sc A.~Taghvaei, J.~de~Wiljes, P.~G. Mehta, and S.~Reich}, {\em {Kalman filter
  and its modern extensions for the continuous-time nonlinear filtering
  problem}}, Journal of Dynamic Systems, Measurement and Control, 140 (2017),
  \url{https://doi.org/10.1115/1.4037780}.

\bibitem{Taghvaei2019}
{\sc A.~Taghvaei, P.~G. Mehta, and S.~P. Meyn}, {\em {Diffusion map-based
  algorithm for gain function approximation in the feedback particle filter}},
  SIAM/ASA Journal on Uncertainty Quantification, 8 (2020), pp.~1090--1117,
  \url{https://doi.org/10.1137/19M124513X}.

\bibitem{sr:tippett03}
{\sc M.~K. Tippett, J.~L. Anderson, C.~H. Bishop, T.~M. Hamill, and J.~S.
  Whitaker}, {\em {Ensemble square root filters}}, Monthly Weather Review, 131
  (2003), pp.~1485--1490,
  \url{https://doi.org/10.1175/1520-0493(2003)131<1485:ESRF>2.0.CO;2}.

\bibitem{Tong2016}
{\sc X.~T. Tong, A.~J. Majda, and D.~Kelly}, {\em {Nonlinear stability of the
  ensemble Kalman filter with adaptive covariance inflation}}, Communications
  in Mathematical Sciences, 14 (2016), pp.~1283--1313,
  \url{https://doi.org/10.4310/CMS.2016.v14.n5.a5}.

\bibitem{VanHandelthesis2007}
{\sc R.~van Handel}, {\em {Filtering, stability, and robustness}}, {P}h{D}
  {T}hesis, California Institute of Technology, 2007,
  \url{https://thesis.library.caltech.edu/4971/}.

\bibitem{Whiteley2018}
{\sc N.~Whiteley}, {\em {Dimension-free Wasserstein contraction of nonlinear
  filters}},  (2018), pp.~1--13, \url{https://arxiv.org/abs/1708.01582}.

\bibitem{Yang2014}
{\sc T.~Yang, H.~A.~P. Blom, and P.~G. Mehta}, {\em {The continuous-discrete
  time feedback particle filter}}, 2014 American Control Conference,  (2014),
  pp.~648--653, \url{https://doi.org/10.1109/ACC.2014.6859259}.

\bibitem{Yang2011}
{\sc T.~Yang, P.~G. Mehta, and S.~P. Meyn}, {\em {Feedback particle filter with
  mean-field coupling}}, Proceedings of the IEEE Conference on Decision and
  Control,  (2011), pp.~7909--7916,
  \url{https://doi.org/10.1109/CDC.2011.6160950}.

\bibitem{Yang2013}
{\sc T.~Yang, P.~G. Mehta, and S.~P. Meyn}, {\em {Feedback particle filter}},
  IEEE Transactions on Automatic Control, 58 (2013), pp.~2465--2480,
  \url{https://doi.org/10.1109/TAC.2013.2258825}.

\bibitem{8100938}
{\sc C.~Zhang, A.~Taghvaei, and P.~G. Mehta}, {\em {Feedback particle filter on
  Riemannian manifolds and matrix Lie groups}}, IEEE Transactions on Automatic
  Control, 63 (2018), pp.~2465--2480,
  \url{https://doi.org/10.1109/TAC.2017.2771336}.

\end{thebibliography}
